\documentclass[11pt]{article}
\usepackage[T1]{fontenc}
\usepackage{lmodern}
\usepackage{amsmath,amssymb,amsfonts,amsthm,mathtools}
\usepackage{mathrsfs,bm,cases,empheq}

\usepackage{indentfirst} 
\usepackage{xcolor}      
\usepackage{graphicx,graphics}
\usepackage{etoolbox,xstring,mfirstuc,textcase}

\usepackage{cite} 

\usepackage[colorlinks=true,linkcolor=blue,citecolor=blue,urlcolor=blue]{hyperref}

\usepackage[capitalize,nameinlink]{cleveref}

\usepackage{changes}
\usepackage{lineno}
\definechangesauthor[name={R1 cusse}, color=orange]{r1}
\definechangesauthor[name={Ty cusse}, color=purple]{ty}

\allowdisplaybreaks[4] 

\newtheorem{theorem}{\textbf{Theorem}}[section]
\newtheorem{lemma}{\textbf{Lemma}}[section]
\newtheorem{proposition}{\textbf{Proposition}}[section]
\newtheorem{corollary}{\textbf{Corollary}}[section]
\newtheorem{remark}{\textbf{Remark}}[section]
\newtheorem{definition}{\textbf{Definition}}[section]
\newtheorem{assumption}{\textbf{Assumption}}[section]

\numberwithin{equation}{section}

\def\be{\begin{equation}}
\def\ee{\end{equation}}
\def\bea{\begin{eqnarray}}
\def\eea{\end{eqnarray}}
\def\bt{\begin{theorem}}
\def\et{\end{theorem}}
\def\bl{\begin{lemma}}
\def\el{\end{lemma}}
\def\bnum{\begin{numcases}{}}
\def\enum{\end{numcases}}
\def\br{\begin{remark}}
\def\er{\end{remark}}
\def\bp{\begin{proposition}}
\def\ep{\end{proposition}}
\def\bc{\begin{corollary}}
\def\ec{\end{corollary}}
\def\bd{\begin{definition}}
\def\ed{\end{definition}}
\def\ba{\begin{assumption}}
\def\ea{\end{assumption}}

\newcommand{\InnerNormal}[2]{( #1, #2 )}
\newcommand{\InnerBig}[2]{\big( #1, #2 \big)}
\newcommand{\InnerMax}[2]{\Big( #1, #2 \Big)}

\newcommand{\PairingNormal}[2]{\langle #1,#2 \rangle}
\newcommand{\PairingBig}[2]{\big\langle #1,#2 \big\rangle}

\newcommand{\NormNormal}[1]{\| #1 \|}
\newcommand{\NormBig}[1]{\big\| #1 \big\|}

\newcommand{\BracketNormal}[1]{( #1 )}
\newcommand{\BracketBig}[1]{\big( #1 \big)}
\newcommand{\BracketMax}[1]{\Big( #1 \Big)}

\newcommand{\AbsNormal}[1]{| #1 |}
\newcommand{\AbsBig}[1]{\big| #1 \big|}
\newcommand{\AbsMax}[1]{\Big| #1 \Big|}

\newcommand{\SetNormal}[1]{\{ #1 \}}
\newcommand{\SetBig}[1]{\big\{ #1 \big\}}
\newcommand{\SetMax}[1]{\Big\{ #1 \Big\}}

\def\sp{\bm{u}}

\begin{document}
\title{ The Cahn--Hilliard--Darcy--Forchheimer system with surfactant: \\ Existence and long-time behavior of global weak solutions}

\author{
Maurizio Grasselli
\thanks{Dipartimento di Matematica,
Politecnico di Milano,
20133 Milano, Italy. \texttt{maurizio.grasselli@polimi.it}.}
\and
Bohan Ouyang
\thanks{School of Mathematical Sciences,
Fudan University,
200433 Shanghai, P.R. China. \texttt{22110180033@m.fudan.edu.cn}.}
\and
Andrea Poiatti
\thanks{Dipartimento di Scienze Matematiche, Fisiche e Informatiche, Università degli Studi di Parma, 43124 Parma, Italy,  \texttt{andrea.poiatti@unipr.it}.}
\and
Hao Wu
\thanks{School of Mathematical Sciences, Fudan University,
200433 Shanghai, P.R. China.  \texttt{haowufd@fudan.edu.cn}.}
}

\date{\today}

\maketitle
\begin{abstract}
\noindent
We consider a diffuse-interface model for two-phase incompressible viscous flows with a soluble surfactant in a bounded porous medium. This hydrodynamic system consists of a Darcy--Forchheimer equation for the seepage velocity $\sp$ coupled with two Cahn--Hilliard equations involving Flory--Huggins type singular potentials, one for the phase-field variable $\phi$, the difference in volume fractions of the two fluids, and the other for the surfactant concentration $\psi$. We study the initial boundary value problem in two or three dimensions, with impermeability boundary conditions for $\sp$ and homogeneous Neumann boundary conditions for $(\phi, \psi)$ and their associated chemical potentials. First, we establish the existence of global weak solutions via an implicit-explicit time-discretization scheme based on the energy dissipation law. Furthermore, applying the seminal results of the first and third authors (arXiv:2510.17296), we prove that every weak solution satisfying an energy inequality converges to a single equilibrium as time tends to infinity. In sharp contrast with the available literature on similar models, in this case weak solutions are enough to guarantee the uniqueness of asymptotic limits, without the necessity of any further eventual regularization.


\end{abstract}

\section{Introduction} \label{intro}
\setcounter{equation}{0}

Multi-phase flows in porous media arise in a wide range of applications, for instance,  enhanced oil recovery, groundwater remediation, filtration processes, and the transport of biological or chemical agents in heterogeneous substrates (see, e.g., \cite{DJ2017,J2023,RK2012} and references therein). In these contexts, the presence of surfactants is crucial. By reducing interfacial tension and modifying wettability, surfactants strongly influence phase separation, interface dynamics, and macroscopic flow properties. A mathematically consistent description of such phenomena therefore requires the coupling of hydrodynamics, phase-field models, and interfacial interactions (see, e.g., \cite{EDAT2013,KK1997,La1992,LK2012}).

Diffuse-interface approaches have become a popular and efficient framework for modeling multi-phase systems with complex interfacial behavior \cite{AMW1998,DF2020}. In this setting, sharp interfaces between mixture components are replaced by thin transition layers in which physical variables vary continuously but steeply. By working with an appropriate order parameter (the phase-field variable), the diffuse-interface model reduces the problem to a set of partial differential equations posed on the entire domain. This allows topological changes in the fluid mixture while avoiding explicit interface tracking. For two-phase flows, the phase-field variable is commonly governed by the Cahn--Hilliard equation, which arises as a (mass-conserving) gradient flow of an underlying free-energy functional (see, e.g., \cite{M2019,Wu2022} and references therein). When surfactants are considered, additional order parameters and coupling terms must be introduced to capture adsorption, diffusion, and interactions with fluid components.
If the multi-phase flow takes place in a porous medium, the momentum balance is typically described by Darcy’s law \cite{GGW2018,G2020}. However, at higher velocities or in more complex porous structures, inertial effects become non-negligible, and Darcy’s law must be augmented by appropriate corrections. A typical extension is given by the Darcy--Forchheimer law, which incorporates a superlinear drag term (see \cite{MTT2016,TT2024,Wh1996} and references therein).
With a similar idea, Cahn--Hilliard--Darcy--Forchheimer systems have recently been applied in tumor-growth modeling, see \cite{BF2025, FLOW2019}. However, a rigorous analytical treatment of such systems within coupled phase-field frameworks remains limited. In particular, the combined presence of Forchheimer drag, surfactant-induced coupling, and singular free-energy potentials has not yet been systematically addressed. The present work aims to fill this gap.

From a mathematical point of view, the nonlinear Forchheimer term facilitates the higher-order integrability of $\sp$, which helps the analysis; on the other hand, it may also pose significant analytical challenges due to its potentially degenerate structure. As an example, we recall that the uniqueness of global weak solutions to the Navier--Stokes system remains an open problem in three dimensions in general. Nevertheless, for the Navier--Stokes system with an additional Forchheimer correction term (i.e., the so-called convective Brinkman--Forchheimer system), the enhanced integrability of the velocity field yields the uniqueness of Leray--Hopf type weak solutions (see \cite{SZ2025} and references therein, see also \cite{BDM2025}).

In this study, we analyze a diffuse-interface model for a two-phase incompressible viscous flow with a soluble surfactant in a bounded porous medium. Assume that $T>0$ is the final time, $\Omega \subset \mathbb{R}^d$, $d \in \SetNormal{2,3}$, is a bounded domain with a smooth boundary $\partial \Omega$. The hydrodynamic system under investigation consists of a Darcy--Forchheimer equation for the (volume averaged) seepage velocity $\sp:\Omega\times [0,T]\to \mathbb{R}^d$, coupled with two Cahn--Hilliard equations: one governing the phase-field variable $\phi:\Omega\times [0,T]\to [-1,1]$ that represents the difference between the volume fractions of the two fluids, and the other describing the surfactant concentration $\psi:\Omega\times [0,T]\to [0,1]$. More precisely, the system reads as follows
\begin{numcases}{}
    \alpha \partial_t \sp + \nu (\phi,\psi) \sp + \eta (\phi,\psi) |\sp|^{r-2} \sp + \nabla \pi = \mu_{\phi} \nabla \phi + \mu_{\psi} \nabla \psi \label{eq:chdf1}\\
    \mathrm{div} \, \sp =0 \label{eq:chdf2}\\
    \partial_t \phi + \sp \cdot \nabla \phi + \sigma_1(\phi) \BracketBig{\overline{\phi} -c} = \mathrm{div} (m_{\phi} (\phi) \nabla \mu_\phi) \label{eq:chdf3}\\
    \mu_{\phi} = -\Delta \phi + \sigma_2 \mathcal{N} \BracketBig{\phi - \overline{\phi}} + F_{\phi}'(\phi) + \partial_{\phi} G(\phi,\psi) \label{eq:chdf4}\\
    \partial_t \psi + \sp \cdot \nabla \psi = \mathrm{div} (m_{\psi} (\psi) \nabla \mu_\psi) \label{eq:chdf5}\\
    \mu_{\psi} = -\beta \Delta \psi + F_{\psi}'(\psi) + \partial_\psi G(\phi,\psi) \label{eq:chdf6}
\end{numcases}
in $\Omega \times (0,T)$, where $r \in (2,\infty)$ is a given constant. The system \eqref{eq:chdf1}--\eqref{eq:chdf6} is subject to the following boundary and initial conditions:
\begin{numcases}{}
    \sp \cdot \mathbf{n} = \partial_{\mathbf{n}} \phi = \partial_{\mathbf{n}} \mu_{\phi} = \partial_{\mathbf{n}} \psi = \partial_{\mathbf{n}} \mu_{\psi} = 0, & on $\partial \Omega \times (0,T)$,\label{eq:chdfb}\\
    \sp|_{t=0} = \sp_0(x), & in $\Omega$, if $\alpha > 0$, \nonumber \\
    \phi|_{t=0} = \phi_0(x), \ \psi|_{t=0} = \psi_0(x),  & in $\Omega$.\label{eq:chdfi}
\end{numcases}
Here, $\mathbf{n}$ denotes the outward unit normal vector on $\partial \Omega$. The constant $\alpha \ge 0$ is a relaxation coefficient (see, e.g., \cite{HW2018} and the references therein), $\nu$ is the Darcy drag coefficient (representing the ratio of the dynamic viscosity to the effective permeability, hereafter referred to as the Darcy coefficient), $\eta$ is the so-called Forchheimer coefficient (related to the porosity of the material), $m_\phi$ and $m_\psi$ are (non-degenerate) mobilities. The term $\sigma_1(\phi) (\overline{\phi} -c)$ in \eqref{eq:chdf3} accounts for possible chemical reactions, where $\sigma_1$ is a non-negative function, $\overline{\phi}=|\Omega|^{-1}\int_\Omega \phi\,\mathrm{d}x$ denotes the spatial mean value of $\phi$, and $c\in (-1,1)$. Here and in the sequel, $\vert\Omega\vert$ denotes the $d$-Lebesgue measure of $\Omega$.

The complex dynamics of the fluid mixture is characterized by the total energy:
\begin{align}
E_{\mathrm{tot}} (\sp,\phi,\psi) \stackrel{\rm{def}}{=} \frac{\alpha}{2} \NormNormal{\sp}_{L^2(\Omega)}^2 + E_{\mathrm{free}} (\phi,\psi),
\label{E-tot}
\end{align}
where the free energy $E_{\mathrm{free}}$ is defined as
\begin{align}
E_{\mathrm{free}}(\phi,\psi)
& \stackrel{\rm{def}}{=}
\int_{\Omega} \left( \frac12 |\nabla \phi|^2 +F_\phi (\phi) + \frac{\sigma_2}{2}\vert \nabla\mathcal{N}(\phi - \overline{\phi})\vert^2 \right)\,\mathrm{d}x
\notag \\
&\quad\ +\int_{\Omega} \left(\frac{\beta}{2}|\nabla \psi|^2+F_\psi (\psi) \right)\,\mathrm{d}x
+\int_{\Omega}G (\phi,\psi)\,\mathrm{d}x.
\label{E-free}
\end{align}
For the sake of simplicity, we have set the width of the diffuse interface to $1$, since the sharp-interface limit is beyond the scope of this study. The free energy $E_{\mathrm{free}}$ enforces physically relevant constraints on the phase-field variables $\phi, \psi$ (see \cite{OGW2026} for a detailed discussion). The chemical potentials $\mu_\phi, \mu_\psi:\Omega \times [0,T]\to \mathbb{R}$ are defined as variational derivatives of $E_{\mathrm{free}}$ corresponding to $\phi, \psi$, respectively.
In \eqref{E-free}, $\beta>0$ is a regularizing coefficient, and $F_\phi$, $F_\psi$ are Flory--Huggins type potential functions.
The bi-variate function $G$ is a coupling term of polynomial type, representing the adsorption and desorption kinetics of the soluble surfactant at the fluid interface.
In addition, the term $\frac{\sigma_2}{2}\vert \nabla\mathcal{N}(\phi - \overline{\phi})\vert^2$ with $\sigma_2\geq 0$, represents nonlocal interactions between fluid components of Ohta--Kawasaki type, where $\mathcal{N}$ stands for the inverse of the Laplace operator subject to a homogeneous Neumann boundary condition.

The main purpose of this work is twofold:
\begin{itemize}
    \item establish the existence of global weak solutions to problem \eqref{eq:chdf1}--\eqref{eq:chdfi}, see Theorem \ref{ws_chdf_thm};
    \item show that each global weak solution converges to a single equilibrium as $t\to \infty$, see Theorem \ref{ws_chdf_convergence_singleton}.
\end{itemize}
We note that related diffuse-interface models involving surfactants in porous media have mainly been analyzed numerically so far (see, for example, \cite{PHRYY2025,Y2022,WT2025}). Our contribution is a first attempt to develop a solid theoretical framework for these problems.

The subsequent analysis is based on two fundamental properties of the coupled system, namely, the mass and energy balances. Compared with the classical Cahn--Hilliard--Darcy system \cite{GGW2018,G2020}, the Forchheimer drag term makes the passage to the limit in the hydrodynamic equation more delicate. Its nonlinear nature prevents a direct identification of the limit from weak convergence of the velocity field $\sp$. To overcome this difficulty, we exploit the monotonicity structure and apply Minty’s trick.
Another issue concerns the control of the convective term in the Cahn--Hilliard equations. The $L^2_t$-$(H^1(\Omega))'$-bound of this term, which is essential for the energy inequality and time regularity of the phase-field variables $\phi, \psi$, requires the $L^{\frac{2r}{r-2}}_tL^{\frac{2r}{r-2}}_x$-estimates of $\phi, \psi$. This becomes particularly delicate in three dimensions when $r$ approaches $2$. In the case of singular potentials, such bounds follow from the $L^2_tL_x^2$-integrability of $F_\phi'(\phi)$ and $F_\psi'(\psi)$. However, if the singular potentials are approximated by regular ones at the level of the approximate system, this boundedness may fail to be uniform with respect to the approximating parameters, thereby preventing a direct compactness argument. For this reason, we adopt an implicit-explicit time-discretization approach that retains singular potentials. This choice ensures that the approximate phase-field variables take their values in the physically admissible ranges. This allows us to derive uniform estimates at the discrete level.

Concerning the long-time behavior of solutions, the standard arguments in the literature, first introduced in \cite{AW2007} for singular potentials, do not seem applicable, as they would require some kind of regularization properties. There are two main obstructions here.
First, due to the highly nonlinear nature of the non-degenerate (continuous) coefficients, like mobilities and permeabilities, the existence of global strong solutions is not expected. Secondly, a weak-strong uniqueness principle appears to be lacking for solutions of the current model.
Consequently, an eventual uniform strict separation property (cf. \cite{AW2007}), which is essential for the application of the {\L}ojasiewicz--Simon inequality in the case of singular potentials, fails to hold.
In the recent contribution \cite{GP2025}, the authors propose a novel argument to prove that global weak solutions to Cahn--Hilliard equations (and even to some Navier--Stokes--Cahn--Hilliard systems) do converge to a unique stationary state. The argument is based on a careful study of times of hyperdissipation of the energy, called ``bad times'', and of times where the $L^2$-norm of the gradient of the chemical potential is bounded, defined as ``good times''. This distinction allows the application of the {\L}ojasiewicz--Simon inequality only along good times (asymptotically). Combined with the energy hyperdissipation of bad times and the validity of an energy inequality, it is sufficient to conclude that each weak trajectory converges to a single equilibrium. Note that no regularization, and thus no weak-strong uniqueness principle, is required. In this study, we aim at crucially extending the arguments in \cite{GP2025} to the more complex case of global weak solutions to problem \eqref{eq:chdf1}--\eqref{eq:chdfi}. One of the main difficulties arises from the lack of mass conservation, which necessitates a nontrivial adaptation of the arguments. Moreover, the regularity of the velocity $\sp$ changes according to the values of $\alpha\geq 0$, and thus its convergence as $t\to\infty$ is nontrivial.
When $\alpha=0$, we can only obtain strong convergence to zero of the translations $\sp(\cdot+t)$ as $t\to\infty$ in suitable Bochner spaces. This is because, due to the lack of (weak) time continuity, in general, there is no pointwise definition of $\sp(t)$ for all $t\geq0$. On the other hand, when $\alpha>0$, thanks to the improved time regularity of $\sp$, we can obtain that $\sp(t)\to\bm{0}$ strongly in $\bm L^2_\sigma(\Omega)$ as $t\to\infty$.

The remainder of this paper is organized as follows. The main results are stated in Section \ref{mr}, after some preliminaries.
In Section \ref{proof_chdf_ws}, we prove the existence of global weak solutions. Section \ref{proof_convergence_equilibrium} is devoted to the analysis of the long-time behavior of weak solutions.

\section{Main Results}\label{mr}
\setcounter{equation}{0}

\subsection{Preliminaries}\label{pre}
We first introduce the notation and conventions used throughout this paper.
Let $\mathcal{X}$ be a (real) Banach space with the norm $\left\|\cdot\right\|_{\mathcal{X}}$.
We denote by $\mathcal{X}^*$ its dual space and by $\langle \cdot,\cdot \rangle_{\mathcal{X}^*,\mathcal{X}}$ the associated duality pairing.
Given a (real) Hilbert space $\mathcal{H}$, its inner product is denoted by $(\cdot,\cdot)_{\mathcal{H}}$.
We assume that $\Omega \subset \mathbb{R}^d$, $d\in\{2,3\}$, is a bounded domain with a sufficiently smooth boundary $\partial \Omega$.
For the standard Lebesgue and Sobolev spaces in $\Omega$, we use the notation $L^{p}(\Omega)$, $W^{k,p}(\Omega)$ for any $p \in [1,\infty]$ and $k\in \mathbb{N}$,
equipped with the corresponding norms $\|\cdot\|_{L^{p}(\Omega)}$, $\|\cdot\|_{W^{k,p}(\Omega)}$, respectively.
When $p = 2$, we use the convention $H^{k}(\Omega) \stackrel{\rm{def}}{=} W^{k,2}(\Omega)$. For simplicity, the norm and the inner product of $L^2(\Omega)$ will be indicated by $\|\cdot\|$ and $(\cdot,\cdot)$, respectively, while the duality pairing between $H^1(\Omega)$ and $H^1(\Omega)^*$ will be indicated by $\langle \cdot,\cdot \rangle$.
The space $W^{k,p}_0 (\Omega)$ denotes the closure of $C^{\infty}_0 (\Omega)$ in $W^{k,p}(\Omega)$.
Its dual space is denoted by $W^{-k,p}(\Omega) \stackrel{\rm{def}}{=} \BracketBig{ W^{k,p'}_0(\Omega) }^*$, where $\frac{1}{p}+\frac{1}{p'}=1$, $p\in (1,\infty)$.
The $L^2$-Bessel potential spaces are denoted by $H^s (\Omega)$, $s\in \mathbb{R}$, which are defined by the restriction of the distributions in $H^s(\mathbb{R}^d)$ to $\Omega$.
For $s >0$, $H^s_0 (\Omega)$ denotes the closure of $C^{\infty}_0(\Omega)$ in $H^s (\Omega)$.
Bold letters will be used for vector-valued spaces, for example, $\bm{L}^p(\Omega)= L^p(\Omega;\mathbb{R}^d)$, $p \in [1,\infty]$. To avoid ambiguity, we also use the notation $(\cdot,\cdot)_{\mathcal{D}}$ to point out the integration domain $\mathcal{D}$.
%

Given a measurable set $I$ of $\mathbb{R}$, we introduce the function space $L^p(I;\mathcal{X})$ with $p\in [1,\infty]$,
which consists of Bochner measurable $p$-integrable functions (if $p \in [1,\infty)$) or essentially bounded functions (if $p =\infty$) with values in a given Banach space $\mathcal{X}$.
If $I=(a,b)$, we write $L^p(a,b;\mathcal{X})$.
In addition, $f \in L^p_{\mathrm{loc}}\BracketBig{[0,\infty);\mathcal{X}}$ if and only if $f \in L^p \BracketBig{0,T;\mathcal{X}}$ for every $T>0$.
The space $L^p_\mathrm{{uloc}}([0,\infty);\mathcal{X})$ denotes the uniformly local variant of $L^p(0,\infty;\mathcal{X})$ consisting of all strongly measurable $f:[0,\infty)\to \mathcal{X}$ such that
\begin{equation*}
    \|f\|_{L^p_\mathrm{{uloc}}([0,\infty);\mathcal{X})} \stackrel{\rm{def}}{=} \mathop\mathrm{sup}\limits_{t\ge0} \|f\|_{L^p(t,t+1;\mathcal{X})} < \infty.
\end{equation*}
If $T\in(0,\infty)$, we have $L^p_\mathrm{{uloc}}([0,T);\mathcal{X})\stackrel{\rm{def}}{=}L^p(0,T;\mathcal{X})$.
For $k\in \mathbb{N}_+$, $p\in [1,\infty)$, we can similarly define $W^{k,p}\BracketBig{0,T;\mathcal{X}}$ and  $W^{k,p}_{\mathrm{uloc}}\BracketBig{[0,\infty);\mathcal{X}}$. When $p=2$, we set
$H^k (0,T;\mathcal{X})=W^{k,2} (0,T;\mathcal{X})$ and $H^k_{\mathrm{uloc}} ([0,\infty);\mathcal{X}) = W^{k,2}_{\mathrm{uloc}} ([0,\infty);\mathcal{X})$.
Let $I=[0,T]$ if $T \in (0,\infty)$ or $I=[0,\infty)$ if $T=\infty$.
Then $BC (I;\mathcal{X})$ denotes the Banach space of all bounded and continuous functions $f:I \to \mathcal{X}$ equipped with the supremum norm and $BUC(I;\mathcal{X})$ is the subspace of all bounded and uniformly continuous functions.
Furthermore, for every $k \in \mathbb{N_+}$, $BC^k(I;\mathcal{X})$ (resp. $BUC^k(I;\mathcal{X})$) denotes all functions $f \in BC(I;\mathcal{X})$ (resp. $f \in BUC(I;\mathcal{X})$), whose Fr\'echet derivatives of order no larger than $k$ exist and belong to $BC(I;\mathcal{X})$ (resp. $BUC(I;\mathcal{X})$).
Finally, we denote by $BC_{w}(I;\mathcal{X})$ the topological vector space of all bounded and weakly continuous functions $f:I \to \mathcal{X}$.

In the subsequent analysis, the following shorthands will be used
\begin{align*}
    & H \stackrel{\rm{def}}{=} L^2\left(\Omega\right), \quad V \stackrel{\rm{def}}{=} H^1\left(\Omega\right), \quad W \stackrel{\rm{def}}{=}\big\{u\in H^2(\Omega)\ |\ \partial_{\mathbf{n}} u=0\ \text{a.e. on}\ \partial\Omega\big\}.
\end{align*}
As usual, $H$ is identified with its dual, and we have the following continuous, dense, and compact embeddings $W \hookrightarrow V  \hookrightarrow H \hookrightarrow V^*$.
Then, we recall the interpolation inequality
\begin{equation*}
    \left\|f\right\|^2 \le \xi \left\|\nabla f\right\|^2 + C\left(\xi\right) \left\|f\right\|^2_{V^*}, \quad \forall f \in V,
\end{equation*}
where $\xi \in (0,1)$ is arbitrary and $C\left(\xi\right)$ is a positive constant depending only on $\xi$ and $\Omega$.

For every $f\in V^*$, $\overline{f}$ denotes its generalized mean value over $\Omega$, given by
$\overline{f}=|\Omega|^{-1}\langle f,1\rangle$ and,
if $f\in L^1(\Omega)$, then we have $\overline{f}=|\Omega|^{-1}\int_\Omega f \, \mathrm{d}x$. Setting
$L^p_{(0)}(\Omega) \stackrel{\rm{def}}{=} \{f\in L^p(\Omega) \ |\ \overline{f} =0\}$, $W^{1,p}_{(0)}(\Omega) \stackrel{\rm{def}}{=} W^{1,p}(\Omega) \cap L^p_{(0)}(\Omega)$, $p \in [1,\infty]$, we recall the well-known Poincar\'{e}--Wirtinger inequality:
\begin{equation*}
\left\| f \right\|_{L^p(\Omega)}\leq C_p \|\nabla f\|_{L^p(\Omega)},\quad \forall\,
f\in W^{1,p}_{(0)}(\Omega),
\end{equation*}
where the positive constant $C_p$ depends only on $p$ and $\Omega$.

Let $A_N: V \to V^*$ denote the extension of the Laplace operator $-\Delta$ subject to the homogeneous Neumann boundary condition, namely,
$\left\langle A_N f,g \right\rangle = \int_{\Omega} \nabla f \cdot \nabla g \, \mathrm{d}x$,  for any $f,g \in V$.
For any given $a \in \mathbb{R}$, we set
\begin{align*}
    & H_{(a)}\stackrel{\rm{def}}{=}\{f\in H\ |\ \overline{f} =a\},\quad
    V_{(a)} \stackrel{\rm{def}}{=} V\cap H_{(a)},\quad
    W_{(a)} \stackrel{\rm{def}}{=} W \cap H_{(a)},\quad
    V^*_{(a)} \stackrel{\rm{def}}{=} \left\{ L \in V^* \ \Big| \   \overline{L} = a \right\}.
\end{align*}
Then, the restriction of $A_N$ to $V_{(0)}$ is a linear isomorphism between $V_{(0)}$ and $V^*_{(0)}$.
Thus, we can define the inverse operator $\mathcal{N} \stackrel{\rm{def}}{=} \big(A_N|_{V_{(0)}}\big)^{-1}:V_{(0)}^* \to V_{(0)}$.
Define
\begin{align*}
    &\left\|L\right\|_* \stackrel{\rm{def}}{=} \left\| \nabla \left(\mathcal{N}L\right)\right\| = \sqrt{\left\langle L,\mathcal{N}L \right\rangle},\quad \forall\, L\in V^*_{(0)}, \quad
    \text{and} \quad
    \left\|L\right\|_{-1}^2 \stackrel{\rm{def}}{=} \left\|L-\overline{L}\right\|_*^2 + |\overline{L}|^2, \quad \forall\, L\in V^*.
\end{align*}
Then $\left\|\cdot\right\|_*$ and $\left\|\cdot\right\|_{-1}$ are equivalent norms in $V^*_{(0)}$ and $V^*$ with respect to the usual dual norms.

Consider $C^{\infty}_{0,\sigma}(\Omega) \stackrel{\rm{def}}{=} \{ \bm{f} \in C^{\infty}_0(\Omega)^d \ | \ \mathrm{div} \bm{f} =0 \}$. We define the following spaces
\begin{align*}
    &\bm{L}_{\sigma}^{p}(\Omega) \stackrel{\rm{def}}{=} \overline{C^{\infty}_{0,\sigma}(\Omega)}^{\bm L^p(\Omega)},
    \quad
    \bm{H}^s_{0,\sigma}(\Omega) \stackrel{\rm{def}}{=} \overline{C^{\infty}_{0,\sigma}(\Omega)}^{\bm H^s(\Omega)},
    \quad
    \bm{W}^{k,p}_{0,\sigma}(\Omega) \stackrel{\rm{def}}{=} \overline{C^{\infty}_{0,\sigma}(\Omega)}^{\bm W^{k,p}(\Omega)},
\end{align*}
for all $p \in (1,\infty)$, $k \in \mathbb{N}_+$ and $s > 0$.
The Helmholtz projection corresponding to $\bm{L}^2_{\sigma}(\Omega)$ is denoted by $\mathbb{P}_{\sigma}$ (cf. \cite{S2001}). We note that $\mathbb{P}_{\sigma} \bm{f} = \bm{f} - \nabla p$, where $p \in V_{(0)}$ is the solution to the (weak) Neumann problem
$\InnerNormal{\nabla p}{\nabla \phi}_{\Omega}
= \InnerNormal{\bm{f}}{\nabla \phi}_{\Omega}$ for any $\phi \in C^{\infty} \BracketBig{\overline{\Omega}}$.
In addition, if $\partial \Omega$ is sufficiently smooth (for example, of class $C^2$), then we can also define the Helmholtz--Weyl decomposition of $L^p(\Omega)$ for all $p \in (1,\infty)$ thanks to the well-posedness of the corresponding (homogeneous) Neumann boundary value problem (see \cite[Section III.1]{G2011}).
When there is no confusion, $\mathbb{P}_{\sigma}$ also denotes the unique projection operator from $\bm{L}^p(\Omega)$ to $\bm{L}_{\sigma}^{p}(\Omega)$, whose null space is
$ \{ \bm{\omega} \in \bm{L}^p(\Omega)\ |\  \bm{\omega} = \nabla p, \ \mathrm{for} \ \mathrm{some} \ p \in W^{1,p}(\Omega) \}$.

The capital letter $C$ denotes a generic positive constant that depends on the structural data of the problem. Its meaning may change from line to line and even within the same chain of computations. Specific dependence will be pointed out if necessary.
Finally, for simplicity, we set
$Q_{(s,t)}=\Omega \times (s,t)$, $Q_{t}=Q_{(0,t)}$, $Q=Q_{(0,\infty)}$, and analogously $S_{(s,t)} = \partial \Omega \times (s,t)$, $S_t = S_{(0,t)}$, $S=S_{(0,\infty)}$.

\subsection{Statement of results}
\label{sectionassumptions}
Let us introduce the assumptions that will be used throughout this paper (cf. \cite{ADG2013ndm,AGG2012,B1999}).

\ba \label{assumptions} \rm
We assume that $\Omega \subset \mathbb{R}^{d}\ (d\in \{ 2,3\})$ is a bounded domain with a smooth boundary $\partial\Omega$ (here it suffices to assume that $\partial \Omega$ is of class $C^2$) and $\beta>0$, $\sigma_2\geq 0$, $c\in (-1,1)$, $\alpha \ge 0$ are given constants.
In addition, we impose the following conditions:
\begin{itemize}
    \item[(\textbf{H1})] The potential function
    $F_{\phi}\in C([-1,1]) \cap C^2((-1,1))$ satisfies
    \begin{equation*}
        \lim\limits_{s \to -1^+} F_{\phi}'(s) = -\infty,\quad
        \lim\limits_{s \to 1^-} F_{\phi}'(s) = + \infty, \quad
        F_{\phi}''(s) \ge \theta_{\phi} > 0, \quad \forall\, s \in (-1,1),
    \end{equation*}
    where $\theta_\phi$ is a given constant.
    Similarly, $F_{\psi}\in C([0,1]) \cap C^2((0,1))$ satisfies
    \begin{equation*}
        \lim\limits_{s \to 0^+} F_{\psi}'(s) = -\infty,\quad
        \lim\limits_{s \to 1^-} F_{\psi}'(s) = + \infty, \quad
        F_{\psi}''(s) \ge \theta_{\psi} > 0, \quad \forall\, s \in (0,1).
    \end{equation*}
    Besides, we set $F_{\psi}(s) = +\infty $ whenever $|s| > 1$ and $F_{\psi}(s) = +\infty $ whenever $s \in (-\infty,0) \cup (1,+\infty)$. Without loss of generality, we take
     \begin{equation*}
         F_{\phi}(0)= F_{\phi}'(0) = F_{\psi}\left(\frac12\right) = F_{\psi}'\left(\frac12\right) = 0.
     \end{equation*}
    Concerning the bi-variate function $G$, we assume that $G(\cdot,\cdot)\in BUC^2(\mathbb{R}^{2})$.
    \item[(\textbf{H2})] We assume $m_{\phi}(\cdot),m_{\psi}(\cdot) \in BUC(\mathbb{R})$.
    They are bounded from below by
    \[
    \underline{m_{\phi}} \stackrel{\rm{def}}{=} \inf\limits_{s \in \mathbb{R}} \{ m_{\phi}(s) \} >0 \quad  \mathrm{and} \quad \underline{m_{\psi}} \stackrel{\rm{def}}{=} \inf\limits_{s \in \mathbb{R}} \{ m_{\psi}(s) \} >0.
    \]
    Their upper bounds are denoted by $ m_{\phi}^* \stackrel{\rm{def}}{=} \sup\limits_{s \in \mathbb{R}} \{ m_{\phi}(s) \}\quad \text{and}\quad m_{\psi}^* \stackrel{\rm{def}}{=} \sup\limits_{s \in \mathbb{R}} \{ m_{\psi}(s) \}$, respectively.
    Moreover, we assume
     $
     \sigma_1(\cdot) \in BUC(\mathbb{R})$, with $\sigma_1(s) \ge 0$, for all $s\in \mathbb{R}.
     $
    \item[(\textbf{H3})] The Darcy coefficient $\nu$ belongs to $BUC( \mathbb{R}^2 )$ with
    $ \underline{\nu} \stackrel{\rm{def}}{=} \inf\limits_{(s_1,s_2) \in \mathbb{R}^2} \{ \nu(s_1,s_2) \} \in \mathbb{R}_+$.
    Its upper bound is denoted by $\nu^* \stackrel{\rm{def}}{=} \sup\limits_{(s_1,s_2) \in \mathbb{R}^2} \{ \nu(s_1,s_2) \}$.
    \item[(\textbf{H4})] The Forchheimer coefficient $\eta$ satisfies $\eta \in BUC( \mathbb{R}^2 )$ with
    $ \underline{\eta} \stackrel{\rm{def}}{=} \inf\limits_{(s_1,s_2) \in \mathbb{R}^2} \{ \eta(s_1,s_2) \} \in \mathbb{R}_+$.
    Its upper bound is denoted by $\eta^* \stackrel{\rm{def}}{=} \sup\limits_{(s_1,s_2) \in \mathbb{R}^2} \{ \eta(s_1,s_2) \}$.
\end{itemize}
\ea
For the existence of weak solutions to the problem \eqref{eq:chdf1}--\eqref{eq:chdfi}, we need the above assumptions on the mobilities and the coefficients in the Darcy--Forchheimer equation.

\begin{remark} \label{compliance_coefficients&extensions}
\rm
Since $\phi\in [-1,1]$ and $\psi \in [0,1]$, only the values of the coefficients $m_\phi$, $m_\psi$, $\nu$, $\eta$ as well as the coupling energy density function $G$ inside the corresponding physically admissible domains are important. On the other hand, for any given function $f \in C^k(\overline{\mathcal{O}})$, where $\mathcal{O} \subset \mathbb{R}^n$ is a Cartesian product of intervals, we can easily extend it to a function $\widetilde{f} \in BUC^k(\mathbb{R}^n)$, keeping the positivity/non-negativity property. This will be convenient in the subsequent analysis.
\end{remark}

\smallskip
Next, we introduce the notion of global weak solutions.

\bd \label{ws_chdf_def}
Let $T \in (0,\infty]$ and set $I=[0,\infty)$ if $T= \infty$ or $I = [0,T]$ if $T<\infty$.
Suppose that the assumptions (\textbf{H1})--(\textbf{H4}) are satisfied.
Then, for any initial data
$\phi_0,\psi_0 \in V$ such that
$(\phi_0,\psi_0)\in[-1,1]\times[0,1]$ almost everywhere in $\Omega$, $\overline{\phi_0}\in(-1,1)$, $\overline{\psi_0}\in(0,1)$, and $\sp_0\in \bm L^2_\sigma(\Omega)$ if $\alpha>0$, a quintuple $(\sp,\phi,\psi,\mu_\phi,\mu_\psi)$ with the following regularity properties
\begin{align*}
    & \sp \in L^r_{\mathrm{loc}} \BracketBig{ I;\bm{L}^r_{\sigma}(\Omega) } \ \ \mathrm{if}\ \alpha = 0, \\
    & \sp \in {BC \BracketBig{ I;\bm{L}^2_{\sigma}(\Omega) }} \cap L^r_{\mathrm{loc}} \BracketBig{ I;\bm{L}^r_{\sigma}(\Omega) }\cap  W^{1,r'}_{\mathrm{loc}} \BracketBig{ I;\bm{L}^{r'}_{\sigma}(\Omega) } \ \ \mathrm{if} \ \alpha >0,  \\
    & \phi,\psi \in BC \BracketBig{I;V} \cap L^2_{\mathrm{loc}} (I;W) \cap L^{\infty} (Q_T), \quad \partial_t \phi, \partial_t \psi \in L^2_{\mathrm{loc}} (I;V^*), \\
    & \phi \in (-1,1) \ \ \mathrm{and} \ \ \psi \in (0,1) \ \ \text{\rm a.e. in} \ \ Q_T,\\
    & \mu_\phi,\mu_\psi \in L^2_{\mathrm{loc}} \BracketBig{I;V},\quad F_{\phi}'(\phi),F_{\psi}'(\psi) \in L^2_{\mathrm{loc}} (I;H),
\end{align*}
is called a (finite energy) weak solution to problem \eqref{eq:chdf1}--\eqref{eq:chdfi} in $I$, if it satisfies
\begin{align}
    - \InnerNormal{ \alpha \sp }{ \partial_t \bm{\theta} }_{Q_T} + \int_{Q_T} ( \nu (\phi,\psi) \sp + \eta (\phi,\psi) | \sp |^{r-2} \sp ) \cdot \bm{\theta} \ \mathrm{d} (x,\tau) = - \InnerNormal{ \phi \nabla \mu_\phi + \psi \nabla \mu_\psi }{ \bm{\theta} }_{Q_T}, \label{testdf}
\end{align}
for all $\bm{\theta} \in C^{\infty}_0 (Q_T)$ with $\mathrm{div}\, \bm{\theta} = 0$, and
\begin{align}\label{testchphi}
   & \InnerBig{\phi}{\partial_t \zeta}_{Q_T} + \InnerBig{\sp \phi}{\nabla \zeta}_{Q_T} - \int_{Q_T} \BracketBig{ \overline{\phi} -c }  \sigma_1 (\phi) \zeta \ \mathrm{d} (x,\tau) = \InnerBig{m_\phi (\phi) \nabla \mu_\phi}{\nabla \zeta}_{Q_T} ,\\
\label{testchpsi}
    & \InnerBig{\psi}{\partial_t \zeta}_{Q_T} + \InnerBig{\sp \psi}{\nabla \zeta}_{Q_T} = \InnerBig{m_\psi (\psi) \nabla \mu_\psi}{\nabla \zeta}_{Q_T},
\end{align}
for all $\zeta \in C_0^{\infty}\big((0,T);C^1\big(\overline{\Omega}\big)\big)$,
equations \eqref{eq:chdf4} and \eqref{eq:chdf6} are satisfied almost everywhere in $Q_T$, the initial conditions
$(\phi,\psi)|_{t=0} =(\phi_0,\psi_0)$ (and $\sp|_{t=0} = \sp_0$ if $\alpha>0$) are satisfied almost everywhere in $\Omega$.
Moreover, we have the following mass relations

\begin{equation} \label{eq:spatial_average}
    \overline{\phi(t)} = c + (\overline{\phi_0}-c) \exp\Big( -\int_0^t \overline{\sigma_1(\phi)} \ \mathrm{d} \tau \Big) \in (-1,1) \quad \text{and}\quad \overline{\psi(t)} = \overline{\psi_0}, \qquad  \forall\, t \in I,
\end{equation}
and the energy identity
\begin{align}
    & E_{\text{tot}} (\sp(t),\phi(t),\psi(t)) + \int_{Q_t} \nu (\phi,\psi) |\sp|^2 + \eta (\phi,\psi) |\sp|^r + m_{\phi}(\phi) |\nabla \mu_\phi|^2 +m_{\psi}(\psi) |\nabla \mu_\psi|^2 \ \mathrm{d} (x,\tau) \nonumber\\
    & \quad + \int_{Q_t} \sigma_1(\phi) \BracketBig{ \overline{\phi}-c } \mu_\phi \ \mathrm{d} (x,\tau) = E_{\text{tot}} (\sp_0,\phi_0,\psi_0), \quad \forall\, t\in I.\label{energyid}
\end{align}

\ed

Our first result reads as follows.
\bt \label{ws_chdf_thm}
Let Assumption \ref{assumptions} be satisfied.

(1) For any initial data $\psi_0, \phi_0 \in V$ with $(\phi_0,\psi_0) \in [-1,1]\times[0,1]$ almost everywhere in $\Omega$ and $\overline{\phi_0} \in (-1,1)$ and $\overline{\psi_0} \in (0,1)$ (and $\sp_0 \in \bm{L}^2_{\sigma} (\Omega)$ if $\alpha >0$), problem \eqref{eq:chdf1}--\eqref{eq:chdfi} admits a global weak solution $(\sp,\phi,\psi,\mu_\phi,\mu_\psi)$ in $[0,\infty)$ in the sense of Definition \ref{ws_chdf_def}. Moreover, there exists a unique function $\pi \in L^{r'}_{\mathrm{loc}} \BracketBig{0,\infty;W^{1,r'}_{(0)}(\Omega) }$  such that
\[
\alpha \partial_t \sp + \nu \BracketBig{ \phi,\psi } \widetilde{\sp} + \eta \BracketBig{ \phi,\psi} |\sp|^{r-2} \sp + \nabla \pi = \mu_\phi \nabla \phi + \mu_\psi \nabla \psi,\quad \text{a.e. in}\ Q.
\]

(2) For every global weak solution $\BracketBig{ \sp,\phi,\psi,\mu_\phi,\mu_\psi}$ in $I=[0,\infty)$ in the sense of Definition \ref{ws_chdf_def}, if in addition,

\begin{equation} \label{additional_requirements_uniform_boundedness}
    \sigma_1 \BracketBig{\phi} \BracketBig{\overline{\phi}-c} \in L^1\BracketBig{I;L^{\frac{6}{5}}(\Omega)},
\end{equation}
then the solution satisfies
\begin{align*}
    & \sp \in L^2 \BracketBig{ I;\bm{L}^2_{\sigma}(\Omega) } \cap L^r \BracketBig{ I;\bm{L}^r_{\sigma}(\Omega) } \ \ \mathrm{if}\ \alpha = 0,\\
    & \sp\in BC \BracketBig{ I;\bm{L}^2_{\sigma}(\Omega) } \cap L^r \BracketBig{ I;\bm{L}^r_{\sigma}(\Omega) }\cap W^{1,r'} \BracketBig{ I;\bm{L}^{r'}_{\sigma}(\Omega) } \ \ \mathrm{if} \ \alpha >0, \\
    & \phi, \psi \in BC \BracketBig{I;V} \cap L^2_{\mathrm{uloc}} (I;W), \\
    & \mu_\phi, \mu_\psi \in L^2_{\mathrm{uloc}} \BracketBig{I;V}, \ \  \nabla \mu_\phi, \nabla \mu_\psi \in L^2 \BracketBig{I;\bm{L}^2(\Omega)},\\
    & F_{\phi}' \BracketBig{\phi},F_{\psi}' \BracketBig{\psi} \in L^2_{\mathrm{uloc}} (I;H).
\end{align*}
\et
\medskip

Next, we study the long-time behavior of global weak solutions established in Theorem \ref{ws_chdf_thm}. To this end, let us  introduce the notion of equilibrium for the evolution system \eqref{eq:chdf1}--\eqref{eq:chdf6}.
As in \cite[Theorem 4.3]{AW2007}, for the energy functional $\widetilde{E}_\phi:H\to \mathbb{R}$ with the effective domain
\[
\mathcal{D}(\widetilde{E}_\phi)=\big\{\phi\in V\ |\  \phi\in [-1,1]\ \text{a.e. in}\ \Omega\big\}
\]
such that
\[
\widetilde{E}_\phi(\phi)=
\begin{cases}
\displaystyle{\int_\Omega\left(\dfrac12|\nabla \phi|^2+  F_{\phi}(\phi)\right)\ \mathrm{d}x},
&\quad \text{if}\ \phi\in \mathcal{D}(\widetilde{E}_\phi),\\
+\infty,&\quad \text{else},
\end{cases}
\]
we define the domain of its
subgradient as
\[
\mathcal{D}(\partial \widetilde{E}_\phi)=\big\{ \phi\in W\ |\
F_{\phi}'(\phi)\in H,\ F_{\phi}''(\phi)|\nabla \phi|^2\in L^1(\Omega)\big\}.
\]
Similarly, for the energy functional
\[
\widetilde{E}_\psi(\psi)=
\begin{cases}
\displaystyle{\int_\Omega\left(\dfrac{\beta}{2}|\nabla \psi|^2+  F_{\psi}(\psi)\right)\ \mathrm{d}x},
&\quad \text{if}\ \psi\in \mathcal{D}(\widetilde{E}_\psi),\\
+\infty,&\quad \text{else},
\end{cases}
\]
with
\[
\mathcal{D}(\widetilde{E}_\psi)=\big\{\psi\in V\ |\  \psi\in [0,1]\ \text{a.e. in}\ \Omega\big\},
\]
we define
\[
\mathcal{D}(\partial \widetilde{E}_\psi)=\big\{ \psi\in W\ |\
F_{\psi}'(\psi)\in H,\ F_{\psi}''(\psi)|\nabla \psi|^2\in L^1(\Omega)\big\}.
\]
According to \cite[Lemma 4.1]{AW2007}, both $E_\phi$ and $E_\psi$ are proper, lower semi-continuous, convex functionals, with
\[
\partial \widetilde{E}_\phi(\phi) =  \SetNormal{-\Delta \phi + F_\phi' (\phi)}, \quad  \partial \widetilde{E}_\psi(\psi) = \SetNormal{- \beta \Delta \psi + F_\psi' (\psi)},
\]
for all $\phi \in \mathcal{D} (\partial \widetilde{E}_\phi)$ and $\psi \in \mathcal{D} (\partial \widetilde{E}_\psi)$.

\bd \label{chdf_equilibrium_def}
Let $(\phi_\infty,\psi_\infty,\mu_{\phi,\infty},\mu_{\psi,\infty}) \in \mathcal{D} (\partial \widetilde{E}_\phi) \times \mathcal{D} (\partial \widetilde{E}_\psi) \times \mathbb{R}^2$
be a solution to the elliptic problem
\begin{numcases}{}
    - \Delta \phi_\infty + \sigma_2 \mathcal{N}\BracketBig{ \phi_\infty - \overline{\phi_\infty} } + F_\phi'(\phi_\infty) + \partial_\phi G(\phi_\infty,\psi_\infty) = \mu_{\phi,\infty}, & a.e. in $\Omega$, \label{ch_stationary_phi} \\
    - \beta \Delta \psi_\infty + F_\psi'(\psi_\infty) + \partial_\psi G(\phi_\infty,\psi_\infty) = \mu_{\psi,\infty}, & a.e. in $\Omega$, \label{ch_stationary_psi} \\
    \partial_{\bm{n}} \phi_\infty = \partial_{\bm{n}} \psi_\infty = 0, & a.e. on $\partial \Omega$, \label{ch_stationary_boundary}
\end{numcases}
subject to the constraint $\sigma_1(\phi_\infty) (\overline{\phi_\infty}-c)=0$.
Then
\begin{itemize}
    \item if $\alpha=0$, we call $(\phi_\infty,\psi_\infty)$ an \emph{equilibrium} of system \eqref{eq:chdf1}--\eqref{eq:chdf6};
    \item if $\alpha >0$, we call $(\bm{0},\phi_\infty,\psi_\infty)$ an \emph{equilibrium} of system \eqref{eq:chdf1}--\eqref{eq:chdf6}.
\end{itemize}
For convenience, we define the following set
\begin{align*}
    \bm{\mathcal{S}} \stackrel{\mathrm{def}}{=} \SetNormal{(\widetilde{\phi},\widetilde{\psi}) \in \mathcal{D} (\partial \widetilde{E}_\phi) \times \mathcal{D} (\partial \widetilde{E}_\psi)&: \ \text{there exist}\ \mu_{\phi}, \mu_{\psi} \in \mathbb{R} \text{ such that } \\
    & \quad (\widetilde{\phi},\widetilde{\psi},\mu_\phi,\mu_\psi) \text{ solves problem \eqref{ch_stationary_phi}--\eqref{ch_stationary_boundary}} }.
\end{align*}
\ed

As in \cite{GP2025}, we introduce the notion of $\omega$-limit sets for $\alpha\geq 0$:
\begin{itemize}
    \item if $\alpha=0$,
        \begin{align*}
            \omega(\phi,\psi) \stackrel{\mathrm{def}}{=} \SetNormal{ (\widetilde{\phi},\widetilde{\psi}) \in (V \cap L^\infty(\Omega))^2 &: \ \text{there exist}\ t_n \to \infty \ \text{such that}\ \phi(t_n) \rightharpoonup \widetilde{\phi}\\
            &\quad \ \text{and}\ \psi(t_n) \rightharpoonup \widetilde{\psi}\  \text{weakly in }V};
        \end{align*}
    \item if $\alpha >0$,
        \begin{align*}
            \omega(\sp,\phi,\psi) \stackrel{\mathrm{def}}{=}
            & \SetNormal{ (\widetilde\sp,\widetilde{\phi},\widetilde{\psi}) \in \mathbf L^2_\sigma(\Omega) \times (V \cap L^\infty(\Omega))^2  : \ \text{there exist}\  t_n \to \infty \ \text{such that}\\
            &\qquad \  \phi(t_n) \rightharpoonup \widetilde{\phi}\  \text{and}\ \psi(t_n) \rightharpoonup \widetilde{\psi}\ \text{weakly in $V$, $\sp(t_n) \rightharpoonup \widetilde{\sp}$ weakly in $\bm{L}^2_\sigma(\Omega)$} }.
        \end{align*}
\end{itemize}
When $\alpha=0$ we do not need to include the velocity $\sp$ in the definition of $\omega$-limit, because in this case the energy does not involve the kinetic energy, and the velocity only appears as a dissipative term.

Since the analyticity of the nonlinear functions $F_\phi, F_\psi, G$ is crucial for the application of the {\L}ojasiewicz--Simon approach (cf. \cite{AW2007,GP2025,O2025}), we assume the following
\begin{itemize}
    \item[(\textbf{H1*})] $F_\phi$ is real analytic in $(-1,1)$, $F_\psi$ is real analytic in $(0,1)$, and $G$ is real analytic in $(-1,1) \times (0,1)$.
\end{itemize}

We are in a position to state the second result.

\bt \label{ws_chdf_convergence_singleton}
Suppose that $(\sp,\phi,\psi,\mu_\phi,\mu_\psi)$ is a global weak solution to problem \eqref{eq:chdf1}--\eqref{eq:chdfi} in the sense of Definition \ref{ws_chdf_def}, which satisfies property \eqref{additional_requirements_uniform_boundedness}. Then

(1) It holds
\begin{align*}
\sp(t+\cdot) \to  \bm{0} \quad \text{strongly in } L^r\BracketBig{0,T;\bm{L}^r(\Omega)}\ \text{ as }t\to\infty,
\end{align*}
for any $T>0$, and
\begin{align*}
    \lim_{t\to\infty}\NormBig{\sp(t)}=0,\quad\text{ if }\alpha>0.
\end{align*}

(2) It holds $\omega(\phi,\psi) \subset \bm{\mathcal{S}}$, if $\alpha=0$, and  $\omega(\sp,\phi,\psi) \subset \SetNormal{\bm{0}} \times \bm{\mathcal{S}}$, if $\alpha>0$.
The $\omega$-limit set $\omega(\phi,\psi)$ (resp. $\omega(\sp,\phi,\psi)$) is compact in $V \times V$ (resp. $\SetNormal{\bm{0}} \times V \times V$), bounded in $W \times W$ (resp. $\SetNormal{\bm{0}} \times W \times W$), and there exists a small constant $\delta_1 >0$ such that
\begin{equation}\label{separation_equilibrium}
    \NormNormal{\phi_\infty}_{L^\infty(\Omega)} \le 1 - 2\delta_1,\quad \  \NormBig{\psi_\infty - \frac{1}{2}}_{L^\infty(\Omega)} \le \frac{1}{2} - 2\delta_1,
\end{equation}
for all $(\phi_\infty,\psi_\infty) \in \omega(\phi,\psi)$ (resp. $(\bm{0},\phi_\infty,\psi_\infty) \in \omega(\sp,\phi,\psi)$).
Moreover, we have
\begin{align}
    & \lim\limits_{t \to \infty} \mathrm{dist}_{V \times V} \BracketBig{ (\phi(t),\psi(t)),\omega(\phi,\psi) }=0, \quad \text{if $\alpha=0$}; \label{dist_alpha=0} \\
    & \lim\limits_{t \to \infty} \mathrm{dist}_{\bm{L}^2_\sigma(\Omega) \times V \times V} \BracketBig{ (\sp(t),\phi(t),\psi(t)),\omega(\sp,\phi,\psi) }=0, \quad \text{if $\alpha>0$}.\label{dist_alpha>0}
\end{align}

(3) If, in addition, the assumption (\textbf{H1*}) is satisfied and
\begin{equation} \label{spatial_average_decay}
    \sup\limits_{t \ge 0} t^{2(1+\rho)} \int_t^\infty \overline{\sigma_1(\phi)} \ \AbsBig{\overline{\phi}-c} \ \mathrm{d} \tau < \infty , \quad \sup\limits_{t \ge 0} t^{1+\rho} \int_t^\infty \NormNormal{\sigma_1(\phi)}_{L^{\frac{6}{5}}(\Omega)} \AbsBig{\overline{\phi}-c} \ \mathrm{d} \tau < \infty,
\end{equation}
for some constant $\rho >0$,
then there exists a pair $(\phi_\infty,\psi_\infty) \in \bm{\mathcal{S}}$ such that
\[
\omega(\phi,\psi)=\SetNormal{(\phi_\infty,\psi_\infty)}\ \ \text{if}\ \alpha=0,\ \ \text{and}\ \ \omega(\sp,\phi,\psi)=\SetNormal{(\bm{0},\phi_\infty,\psi_\infty)}\ \ \text{if}\ \alpha>0.
\]
In both cases, we have
\[
\lim\limits_{t \to \infty} \big(\NormNormal{ \phi(t) - \phi_\infty }_{H^1(\Omega)} + \NormNormal{ \psi(t) - \psi_\infty }_{H^1(\Omega)}\big) = 0.
\]
\et

\begin{remark}\rm
    For $\alpha \ge 0$, we can show that $\sp \in L^{\frac{r}{2}}(T^*,\infty; \bm{L}^r(\Omega)) \cap L^1(T^*,\infty; \bm{L}^2_\sigma(\Omega))$ for some sufficiently large $T^*>0$ (see \eqref{L1_integrability}).
    When $\alpha=0$, the lack of time regularity for $\sp$ prevents us from concluding its pointwise convergence to $\bm{0}$ in $\bm{L}^2_\sigma(\Omega)$ as $t \to \infty$ directly. Nevertheless, if one can further show that $\sp \in BUC([T^*,\infty); \bm{L}^2_\sigma(\Omega))$, then it follows that
    \begin{align}
        \lim_{t \to \infty} \NormNormal{\sp(t)} = 0. \label{c1}
    \end{align}
\end{remark}
\begin{remark} \label{spatial_average_decay_condition}
\rm
We note that  \eqref{additional_requirements_uniform_boundedness} and \eqref{spatial_average_decay} (cf. \cite[Section 6]{JWZ2015}) are technical assumptions that can be fulfilled, for instance, when $\sigma_1\geq0$ is a constant independent of $\phi$. Indeed, in this case, we have
    \[
    \sigma_1 (\overline{\phi(t)}-c) = \sigma_1 e^{-\sigma_1 t}(\overline{\phi_0}-c),\quad \forall\,t\geq 0,
    \]
    which implies the desired properties of $\sigma_1 (\overline{\phi}-c)$.
    Furthermore, it follows from \eqref{eq:spatial_average} that the requirements \eqref{additional_requirements_uniform_boundedness} and \eqref{spatial_average_decay}  are also satisfied in the following two cases (cf. \cite[Remark 2.4]{OGW2026}):
    either $\overline{\phi_0} = c$, or $\sigma_1= \sigma_0 f_0$, where $\sigma_0$ is a nonnegative constant and $f_0 \in C([-1,1])$ is a positive function.
\end{remark}
\begin{remark}\rm
    The results presented in Theorems \ref{ws_chdf_thm}, \ref{ws_chdf_convergence_singleton} can easily be extended to the
    case $r=2$, i.e., the Darcy's regime, because the uniform $L^\infty$-bounds on the phase-field variables are sufficient to control the convective terms by using  the $L^2$-integrability of the velocity. In this linear drag case, the regularity of the associated pressure $\pi$ can be enhanced to $L^q_{\mathrm{loc}}(0,\infty;W^{1,q}_{(0)}(\Omega))$ for all $q \in (1,2)$, since the restriction imposed by the Forchheimer term in $L^{r'}$ is removed. In addition, the result of Theorem \ref{ws_chdf_convergence_singleton} on the long-time behavior of global weak solutions is new even for the Cahn--Hilliard--Darcy system without surfactant interaction (cf. \cite{G2020,GGW2018}).
\end{remark}

\section{Existence of Global Weak Solutions}
\label{proof_chdf_ws}
\setcounter{equation}{0}

The proof of Theorem \ref{ws_chdf_thm} is based on a suitable implicit-explicit time-discretization scheme in the spirit of  \cite{ADG2013ndm,OGW2026} (see also, e.g., \cite{AbelsGarckePoiatti}). We first construct a family of approximate solutions and then derive uniform estimates for such solutions. After passing to the limit, we can extract a convergent subsequence and establish the existence of a global weak solution.
Although in the case $\alpha = 0$ there is no initial value condition imposed on $\sp$, for convenience of notation, below we set $\sp_0 = \mathbf 0$ if $\alpha = 0$. Nevertheless, this value will never come into play for $\alpha=0$.

\subsection{Implicit-explicit time discretization}\label{implicit_time_discret}
Let us introduce the implicit-explicit time discretization.
In this subsection, we shall work with a slightly stronger assumption on the mobilities, that is, $m_\phi, m_\psi \in BUC^1(\mathbb{R})$. Besides, we set $\widetilde{E}_{\phi,a} = \widetilde{E}_{\phi}|_{L^2_{(a)}(\Omega)}$ and $\widetilde{E}_{\psi,b} = \widetilde{E}_{\psi}|_{L^2_{(b)}(\Omega)}$ for any given $a \in (-1,1)$ and $b \in (0,1)$ (see \cite[Section 3]{OGW2026}).
For any $N \in \mathbb{N}_+$, we define the time step $h=\frac{1}{N}$.
Given
\[
\sp^k \in \bm{L}^2_{\sigma}(\Omega), \ \ \phi^k, \psi^k \in V \  \text{ with } F'_\phi \BracketBig{\phi^k}, F'_\psi \BracketBig{\psi^k} \in H,
\]
we seek functions
\[
(\sp,\phi,\psi,\widehat{ \mu_\phi },\widehat{ \mu_\psi }) = (\sp^{k+1},\phi^{k+1},\psi^{k+1},\widehat{ \mu_\phi }^{k+1},\widehat{ \mu_\psi }^{k+1}),
\]
satisfying
\[
\sp \in \bm{L}^r_{\sigma}(\Omega), \quad
\phi \in \mathcal{D} (\partial \widetilde{E}_{\phi,a}),\quad
\psi \in \mathcal{D} (\partial \widetilde{E}_{\psi,b}),\quad
\widehat{ \mu_\phi }, \widehat{ \mu_\psi } \in W_{(0)},
\]
with prescribed mean values
\begin{align}
a= \overline{\phi^k} - h \overline{\sigma_1(\phi^k)} \BracketBig{\overline{\phi^k}-c}, \quad b =\overline{\psi^k},
\label{coe-ab}
\end{align}
as a solution to the following nonlinear discrete system:
\begin{numcases}{}
    \frac{\alpha}{h} (\sp-\sp^k,\bm{\theta}) + \InnerNormal{ \nu(\phi^k,\psi^k) \sp + \eta (\phi^k,\psi^k) |\sp|^{r-2} \sp }{ \bm{\theta} }\nonumber \\
    \quad = - \InnerNormal{ \phi^k \nabla \widehat{ \mu_\phi } + \psi^k \nabla \widehat{ \mu_\psi } }{ \bm{\theta} },
    & $\forall\, \bm{\theta} \in C^{\infty}_{0,\sigma} (\Omega)$,
    \label{dfdiscrete_project}\\
    \frac{\phi -\phi^k}{h} + \sp \cdot \nabla \phi^k + \sigma_1\BracketBig{\phi^k}
    \BracketBig{\overline{\phi^k} -c} = \mathrm{div} \BracketBig{m_\phi \BracketBig{\phi^k} \nabla \widehat{ \mu_\phi }},
    & a.e. in $\Omega$,
    \label{chdiscrete1_project}\\
    \widehat{ \mu_\phi } - P_0 G_{\phi} \BracketBig{\phi,\phi^k,\psi} - \sigma_2 \mathcal{N} \BracketBig{\phi - \overline{\phi}} = -\Delta \phi + P_0 F_{\phi}'(\phi),
    & a.e. in $\Omega$,
    \label{chdiscrete2_project}\\
    \frac{\psi -\psi^k}{h} + \sp \cdot \nabla \psi^k = \mathrm{div} \BracketBig{m_\psi \BracketBig{\psi^k} \nabla \widehat{ \mu_\psi }},
    & a.e. in $\Omega$,
    \label{chdiscrete3_project}\\
    \widehat{ \mu_\psi } - P_0 G_{\psi} \BracketBig{\phi^k,\psi,\psi^k} = -\beta \Delta \psi + P_0 F_{\psi}'(\psi),
    & a.e. in $\Omega$.
    \label{chdiscrete4_project}
\end{numcases}
As in \cite{OGW2026}, the functions $G_\phi$, $G_\psi$ are defined as follows
\begin{equation}\notag
    G_{\phi}(a,b,c)\stackrel{\rm{def}}{=}\left\{
    \begin{aligned}
        &\frac{G(a,c)-G(b,c)}{a-b} , \ &  \mathrm{if} \ a \neq b,\\
        &\partial_\phi G (a,c), \ &  \mathrm{if} \ a=b,
    \end{aligned}
    \right.
\end{equation}
\begin{equation}\notag
    G_{\psi}(c,a,b)\stackrel{\rm{def}}{=}\left\{
    \begin{aligned}
        &\frac{G(c,a)-G(c,b)}{a-b} , \ &  \mathrm{if} \ a \neq b,\\
        &\partial_\psi G (c,a), \ &  \mathrm{if} \ a=b,
    \end{aligned}
    \right.
\end{equation}
for any $(a,b,c) \in \mathbb{R}^3$.

We first recall the following \textit{a priori} estimates for $(\phi,\mu_\phi)$ obtained from equation \eqref{chdiscrete2_project} (resp. for $(\psi,\mu_\psi)$ from equation \eqref{chdiscrete4_project}) like in \cite{OGW2026,ADG2013ndm}.

\bl \label{apriori_ch}
Assume that $\phi \in \mathcal{D}(\partial \widetilde{E}_\phi)$, $ \psi \in \mathcal{D}(\partial \widetilde{E}_\psi)$ and $\widehat{\mu_\phi},\widehat{\mu_\psi} \in V_{(0)}$ solve equations \eqref{chdiscrete2_project} and  \eqref{chdiscrete4_project} with given functions $\phi^k,\psi^k \in H^2(\Omega)$ satisfying $ \BracketBig{\phi^k,\psi^k} \in [-1,1]\times[0,1]$ in $\Omega$.
For any $\delta \in (0,\tfrac{1}{2})$, there exist two positive constants $K_{\phi} = K_{\phi}(\delta)$, $K_{\psi} = K_{\psi}(\delta)$ such that if
$\big( \overline{\phi},\overline{\psi} \big) \in [-1+\delta,1-\delta] \times [\delta,1-\delta]$, we have
\begin{align*}
    & \NormNormal{F_{\phi}'(\phi)} \le K_{\phi}\BracketNormal{\NormNormal{\widehat{\mu_\phi}} + 1},\qquad \ \ \,    \NormNormal{F_{\psi}'(\psi)} \le K_{\psi}\BracketNormal{\NormNormal{\widehat{\mu_\psi}} + 1},\\
    &\NormNormal{\partial \widetilde{E}_\phi(\phi)} \le K_{\phi} \BracketNormal{\NormNormal{\widehat{\mu_\phi}}+1},\qquad
    \NormNormal{\partial \widetilde{E}_\psi(\psi)} \le K_{\psi} \BracketNormal{\NormNormal{\widehat{\mu_\psi}}+1}.
\end{align*}
\el
\medskip

Let $a,b$ be determined as in \eqref{coe-ab}. We define the following spaces (cf. \cite{OGW2026,ADG2013ndm}):
\begin{align*}
    & X \stackrel{\rm{def}}{=} \mathcal{D} (\partial \widetilde{E}_{\phi,a}) \times \mathcal{D} ( \partial \widetilde{E}_{\psi,b}) \times W_{(0)} \times W_{(0)}, \\
    & Y \stackrel{\rm{def}}{=} H_{(0)} \times H_{(0)} \times H_{(0)} \times H_{(0)}, \\
    & \widetilde{X} \stackrel{\rm{def}}{=} V_{(a)} \times V_{(b)} \times V_{(0)} \times V_{(0)}.
\end{align*}
The following result yields the existence of a solution to the time-discrete system \eqref{dfdiscrete_project}--\eqref{chdiscrete4_project}.

\bp \label{discrete_exist_apriori}
Let $\sp^k \in \bm{L}^2_{\sigma} (\Omega)$ and $\phi^k,\psi^k \in W^{2,6}(\Omega) \cap W$ be given, with $(\phi^k, \psi^k)\in [-1,1] \times [0,1]$ almost everywhere in $\Omega$. Then there exists $(\sp, \phi, \psi, \widehat{ \mu_\phi }, \widehat{ \mu_\psi }) \in \bm{L}^r_\sigma(\Omega) \times X$ that solves  problem \eqref{dfdiscrete_project}--\eqref{chdiscrete4_project} and satisfies the following discrete energy inequality:
\begin{align}
    & E_{\text{tot}} (\sp, \phi, \psi) + \frac{1}{2} \BracketBig{ \NormNormal{ \nabla \phi - \nabla \phi^k }^2 + \beta \NormNormal{ \nabla \psi - \nabla \psi^k }^2 } \nonumber \\
    & \qquad + \frac{\alpha}{2} \NormBig{ \sp - \sp^k }^2 + \frac{\sigma_2}{2} \NormBig{ \phi - \phi^k - \BracketBig{ \overline{\phi} - \overline{\phi^k} } }^2_* \nonumber \\
    & \qquad + h \int_\Omega \BracketBig{ \nu (\phi^k,\psi^k) | \sp |^2 + \eta (\phi^k,\psi^k) | \sp |^r + m_\phi (\phi^k) | \nabla \widehat{ \mu_\phi } |^2 + m_\psi (\psi^k) | \nabla \widehat{ \mu_\psi } |^2 } \ \mathrm{d} x \nonumber \\
    & \qquad + h \BracketBig{ \overline{\phi^k} - c } \int_\Omega \sigma_1(\phi^k) \widehat{ \mu_\phi } \ \mathrm{d} x - |\Omega| \BracketBig{ \overline{\phi}-\overline{\phi^k} } \BracketBig{ \overline{F'_\phi (\phi)} + \overline{G_\phi (\phi,\phi^k,\psi)} } \nonumber \\
    & \quad \le E_{\text{tot}} (\sp^k, \phi^k, \psi^k). \label{discrete_energy_project}
\end{align}
In addition, $\phi, \psi \in W^{2,p}(\Omega)$, where $p=6$ if $d=3$ and $p \in [2,\infty)$ if $d=2$.
\ep

\noindent \textbf{Proof.}
We first derive \textit{a priori} estimates for every solution $(\sp,\phi,\psi,\widehat{ \mu_\phi },\widehat{ \mu_\psi }) \in \bm{L}^r_\sigma(\Omega) \times X$ to problem \eqref{dfdiscrete_project}--\eqref{chdiscrete4_project}.
Recalling the definition of $X$, a density argument justifies taking $\bm{\theta} = \sp$ in \eqref{dfdiscrete_project}.
This gives
\begin{align*}
    & \frac{\alpha}{2} \BracketBig{ \NormNormal{\sp}^2 + \NormNormal{\sp-\sp^k}^2 } + h \int_\Omega \BracketBig{ \nu (\phi^k,\psi^k) | \sp |^2 + \eta (\phi^k,\psi^k) | \sp |^r } \ \mathrm{d} x \\
    & \quad = \frac{\alpha}{2} \NormNormal{ \sp^k }^2 + h \InnerNormal{ \widehat{ \mu_\phi } \nabla \phi^k + \widehat{ \mu_\psi } \nabla \psi^k }{ \sp },
\end{align*}
where we used integration by parts.
On the other hand, a direct computation yields that
\begin{align*}
    & \frac{1}{2} \NormNormal{\nabla \phi}^2 + \frac{\beta}{2} \NormNormal{\nabla \psi}^2 + \frac{\sigma_2}{2} \NormBig{ \phi - \overline{\phi} }^2_* + \int_\Omega \BracketBig{ F_\phi (\phi) + F_\psi (\psi) + G(\phi,\psi) } \ \mathrm{d} x \\
    & \qquad + \frac{1}{2} \NormNormal{\nabla \phi - \nabla \phi^k}^2 + \frac{\beta}{2} \NormNormal{\nabla \psi - \nabla \psi^k}^2 + \frac{\sigma_2}{2} \NormBig{ \phi - \overline{\phi} - \BracketBig{ \phi^k - \overline{\phi^k} } }^2_* \\
    & \qquad + h \int_\Omega \BracketBig{ m_\phi (\phi^k) |\nabla \widehat{ \mu_\phi }|^2 + m_\psi (\psi^k) |\nabla \widehat{ \mu_\psi }|^2 } \ \mathrm{d} x \\
    & \qquad + h \BracketBig{ \overline{\phi^k} - c } \int_\Omega \sigma_1 \BracketBig{\phi^k} \widehat{ \mu_\phi } \ \mathrm{d} x - |\Omega| \BracketBig{ \overline{\phi}-\overline{\phi^k} } \BracketBig{ \overline{F'_\phi (\phi)} + \overline{G_\phi (\phi,\phi^k,\psi)} } \\
    & \qquad + h \InnerNormal{ \widehat{ \mu_\phi } \nabla \phi^k + \widehat{ \mu_\psi } \nabla \psi^k }{ \sp } \\
    & \quad \le \frac{1}{2} \NormNormal{\nabla \phi^k}^2 + \frac{\beta}{2} \NormNormal{\nabla \psi^k}^2 + \frac{\sigma_2}{2} \NormBig{ \phi^k - \overline{\phi^k} }^2_* + \int_\Omega \BracketBig{ F_\phi (\phi^k) + F_\psi (\psi^k) + G(\phi^k,\psi^k) } \ \mathrm{d} x.
\end{align*}
Adding the above two estimates together, we obtain the inequality \eqref{discrete_energy_project}.
Moreover, it follows from \cite[Lemma 7.4]{GGW2018} (see also \cite[Lemma 2]{A2009}) that there exist some positive constants $C_\phi$ and $C_\psi$, which depend on $\Omega$, $\overline{\phi_0}$, $\overline{\psi_0}$, $c$ and $p$ but not on $\phi$, $\psi$, $\widehat{ \mu_\phi }$ and $\widehat{ \mu_\psi }$, such that
\begin{align*}
    & \NormNormal{\phi}_{W^{2,p}(\Omega)} + \NormNormal{F'_\phi(\phi)}_{L^p(\Omega)} \le C_\phi (1 + \NormNormal{\phi} + \NormNormal{\nabla \widehat{ \mu_\phi }}), \\
    & \NormNormal{\psi}_{W^{2,p}(\Omega)} + \NormNormal{F'_\psi(\psi)}_{L^p(\Omega)} \le C_\psi (1 + \NormNormal{\psi} + \NormNormal{\nabla \widehat{ \mu_\psi }}),
\end{align*}
for all $p \in [2,\infty)$ if $d=2$ and $p=6$ if $d=3$.

Next, we prove the existence of a solution to the problem \eqref{dfdiscrete_project}--\eqref{chdiscrete4_project}.
To this end, we split this system according to its convex-concave structure and employ a fixed point argument.
For any $\bm{w} = (\phi,\psi,\widehat{ \mu_\phi },\widehat{ \mu_\psi }) \in X$, we define the mappings $\mathcal{L}_k : X \to Y$ and $\mathcal{F}_k: \widetilde{X} \to Y$ as follows:
\begin{eqnarray*}
    \mathcal{L}_k(\bm{w}) &=&
    \left(
    \begin{array}{c}
           -\mathrm{div} \BracketBig{ m_\phi \BracketBig{\phi^k} \nabla \widehat{ \mu_\phi }} \\
           -\Delta \phi + P_0 F_{\phi}'(\phi) \\
           -\mathrm{div} \BracketBig{ m_\psi \BracketBig{\psi^k} \nabla \widehat{ \mu_\psi }} \\
           - \beta \Delta \psi + P_0 F_{\psi}'(\psi)
    \end{array}
    \right),
\end{eqnarray*}
and
\begin{eqnarray*}
    \mathcal{F}_k(\bm{w})&=&
    \left(
    \begin{array}{c}
        -\dfrac{\phi -\phi^k}{h} - \mathcal{M}^{-1}_k \BracketMax{ \mathcal{C}_k (\widehat{ \mu_\phi },\widehat{ \mu_\psi })  + \dfrac{\alpha}{h} \sp^k } \cdot \nabla \phi^k - \sigma_1 (\phi^k) \BracketBig{\overline{\phi^k} -c} \\
        \widehat{ \mu_\phi } - P_0 G_{\phi} \BracketBig{\phi,\phi^k,\psi} - \sigma_2 \mathcal{N} \BracketBig{ \phi - \overline{\phi} } \\
        -\dfrac{\psi -\psi^k}{h} - \mathcal{M}^{-1}_k \BracketMax{ \mathcal{C}_k (\widehat{ \mu_\phi },\widehat{ \mu_\psi }) + \dfrac{\alpha}{h} \sp^k } \cdot \nabla \psi^k \\
        \widehat{ \mu_\psi } - P_0 G_{\psi} \BracketBig{\phi^k,\psi,\psi^k}
    \end{array}
    \right),
\end{eqnarray*}
where $\mathcal{M}_k: \bm{L}^r_{\sigma} (\Omega) \to \bm{L}^r_{\sigma} (\Omega)^* \cong \bm{L}^{r'}_\sigma (\Omega)$ is a monotone operator defined by
\[
\PairingBig{ \mathcal{M}_k (\sp) }{\bm{v}}_{\bm{L}^r_{\sigma} (\Omega)^*,\bm{L}^r_{\sigma} (\Omega)} = \InnerBig{ \dfrac{\alpha}{h} \sp + \nu (\phi^k,\psi^k) \sp }{\bm{v}} + \PairingNormal{ \eta (\phi^k,\psi^k) |\sp|^{r-2} \sp }{\bm{v}}_{\bm{L}^{r'}(\Omega),\bm{L}^r(\Omega)} ,
\]
for all $\bm{v} \in \bm{L}^r_{\sigma} (\Omega)$, and $\mathcal{C}_k:V \times V \to \bm{L}^2_{\sigma}(\Omega)$ is defined by
\[
\mathcal{C}_k(\mu_\phi,\mu_\psi) = - \mathbb{P}_\sigma ( \phi^k \nabla \mu_\phi + \psi^k \nabla \mu_\psi ).
\]
We observe that $\NormNormal{\mathcal{M}_k (\sp)}_{\bm{L}^{r'}_\sigma (\Omega)} \le C ( \NormNormal{\sp}_{\bm{L}^2_{\sigma} (\Omega)} + \NormNormal{\sp}^{r-1}_{\bm{L}^r_{\sigma} (\Omega)} )$. The injectivity of $\mathcal{M}_k$ follows from its strict monotonicity. In addition, its monotonicity, hemicontinuity and coercivity ensure the surjectivity of $\mathcal{M}_k$ (see \cite[Theorem 9.14-1]{Ciarlet2013}). Moreover, taking $\bm{f}_1,\bm{f}_2 \in \bm{L}^2_{\sigma} (\Omega)$ we have
\begin{align*}
    \NormBig{ \mathcal{M}_k^{-1} \bm{f}_1 - \mathcal{M}_k^{-1} \bm{f}_2 }^2
    & \le \frac{1}{\underline{\nu}} \InnerMax{ \nu (\phi^k,\psi^k) \BracketNormal{ \mathcal{M}_k^{-1} \bm{f}_1 - \mathcal{M}_k^{-1} \bm{f}_2 } }{ \mathcal{M}_k^{-1} \bm{f}_1 - \mathcal{M}_k^{-1} \bm{f}_2 } \\
    & \le \frac{1}{\underline{\nu}} \left\langle \mathcal{M}_k^{-1} \bm{f}_1 - \mathcal{M}_k^{-1} \bm{f}_2, \bm{f}_1 - \bm{f}_2 \right\rangle \\
    & \le \frac{1}{2} \NormNormal{ \mathcal{M}_k^{-1} \bm{f}_1 - \mathcal{M}_k^{-1} \bm{f}_2 }^2 + \frac{1}{2\underline{\nu}^2} \NormNormal{ \bm{f}_1 - \bm{f}_2 }^2,
\end{align*}
which implies the continuity of $\mathcal{M}_k^{-1}|_{\bm{L}^2_{\sigma}(\Omega)}$ from $\bm{L}^2_{\sigma}(\Omega)$ to $\bm{L}^2_{\sigma}(\Omega)$.
Recalling that $\phi^k \in [-1,1]$ almost everywhere in $\Omega$, we have
\[
\AbsMax{ \int_\Omega \phi^k \nabla \mu_\phi \cdot \bm{\theta} \ \mathrm{d} x } \le \NormNormal{\nabla \mu_\phi} \NormNormal{\bm{\theta}},
\]
for all $\bm{\theta} \in C^\infty_{0,\sigma} (\Omega)$, and a similar result holds for $\mu_\psi$.
Consequently, $\mathcal{C}_k$ is continuous and
\[
\NormBig{\mathcal{M}_k^{-1} \BracketBig{ \mathcal{C}_k(\widehat{ \mu_\phi },\widehat{ \mu_\psi }) + \dfrac{\alpha}{h} \sp^k } } \le \dfrac{1}{\underline{\nu}} \BracketBig{ \NormNormal{\nabla \widehat{ \mu_\phi }} + \NormNormal{\nabla \widehat{ \mu_\psi }} + \frac{\alpha}{h} \NormNormal{\sp^k} }.
\]
Therefore, $\mathcal{F}_k$ is continuous from $\widetilde{X}$ to $Y$ and maps bounded sets in $\widetilde{X}$ to bounded sets in $Y$.
On the other hand, arguing as in \cite{OGW2026,ADG2013ndm}, we can show that $\mathcal{L}_k^{-1}$ is a continuous and compact map from $Y$ to $\widetilde{X}$.
In summary, the map $\mathcal{K}_k \stackrel{\rm{def}}{=} \mathcal{F}_k \circ \mathcal{L}^{-1}_k $ is continuous and compact from $Y$ to $Y$.

We see that a tuple $(\sp,\bm{w})$ with $\bm{w} = (\phi,\psi,\widehat{ \mu_\phi },\widehat{ \mu_\psi }) \in X$ solves the problem \eqref{dfdiscrete_project}--\eqref{chdiscrete4_project} if and only if $\bm{f} = \mathcal{L}_k (\bm{w}) \in Y$ is a fixed point of $\mathcal{K}_k$.
Indeed, once a fixed point $\bm{f}$ is obtained, one recovers the corresponding state variables $\bm{w} = \mathcal{L}_k^{-1}(\bm{f})$, which in turn uniquely determine the chemical potentials $\widehat{\mu_\phi}$ and $\widehat{\mu_\psi}$.
Then, the corresponding velocity field is explicitly given by the relation
\[
\sp = \mathcal{M}_k^{-1} \BracketMax{\mathcal{C}_k(\widehat{\mu_\phi},\widehat{\mu_\psi}) + \frac{\alpha}{h} \sp^k }.
\]
By the Leray--Schauder principle, to prove the existence of such a fixed point, it suffices to show that
\[
\text{there exists}\ R >0 \text{ such that, if } \bm{f} \in Y \text{ and } 0 \le \lambda \le 1 \text{ fulfills } \bm{f} = \lambda \mathcal{K}_k \bm{f}, \text{ then } \NormNormal{\bm{f}}_Y \le R.
\]
Let $(\lambda,\bm{f})$ be a solution to $\bm{f} = \lambda \mathcal{K}_k \bm{f}$. Set $\bm{w} = \mathcal{L}_k^{-1}(\bm{f})=(\phi,\psi,\widehat{\mu_\phi},\widehat{\mu_\psi})$. Then we have
\[
\bm{f}=\lambda \mathcal{K}_k \bm{f} \quad \Longleftrightarrow
    \quad \mathcal{L}_k(\bm{w}) - \lambda \mathcal{F}_k(\bm{w})=\mathbf{0},
\]
which is equivalent to the following weak formulation
\begin{numcases}{}
    \frac{\alpha}{h} (\sp-\sp^k,\bm{\theta}) + \InnerNormal{ \nu(\phi^k,\psi^k) \sp + \eta (\phi^k,\psi^k) |\sp|^{r-2} \sp }{ \bm{\theta} }\nonumber \\
    \quad = \InnerNormal{ \widehat{ \mu_\phi } \nabla \phi^k + \widehat{ \mu_\psi } \nabla \psi^k }{ \bm{\theta} },
    & $\forall\, \bm{\theta} \in C^{\infty}_{0,\sigma} (\Omega)$,
    \label{df_schauder}\\
    \lambda \frac{\phi -\phi^k}{h} + \lambda \sp \cdot \nabla \phi^k + \lambda \sigma_1\BracketBig{\phi^k}
    \BracketBig{\overline{\phi^k} -c} = \mathrm{div} \BracketBig{m_\phi \BracketBig{\phi^k} \nabla \widehat{ \mu_\phi }},
    & a.e. in $\Omega$,
    \label{ch1_schauder}\\
    \lambda \widehat{ \mu_\phi } - \lambda P_0 G_{\phi} \BracketBig{\phi,\phi^k,\psi} - \lambda \sigma_2 \mathcal{N} \BracketBig{\phi - \overline{\phi}} = -\Delta \phi + P_0 F_{\phi}'(\phi),
    & a.e. in $\Omega$,
    \label{ch2_schauder}\\
    \lambda \frac{\psi -\psi^k}{h} + \lambda \sp \cdot \nabla \psi^k = \mathrm{div} \BracketBig{m_\psi \BracketBig{\psi^k} \nabla \widehat{ \mu_\psi }},
    & a.e. in $\Omega$,
    \label{ch3_schauder}\\
    \lambda \widehat{ \mu_\psi } - \lambda P_0 G_{\psi} \BracketBig{\phi^k,\psi,\psi^k} = -\beta \Delta \psi + P_0 F_{\psi}'(\psi),
    & a.e. in $\Omega$.
    \label{ch4_schauder}
\end{numcases}
Using a similar argument for \eqref{discrete_energy_project}, we find the following estimate
\begin{align*}
    & h \int_\Omega \left(\lambda \nu\BracketBig{\phi^k,\psi^k}\AbsNormal{\sp}^2 + \lambda \eta \BracketBig{\phi^k,\psi^k}\AbsNormal{\sp}^r + m_\phi\BracketBig{\phi^k} \AbsNormal{\nabla \widehat{\mu_\phi}}^2 + m_\psi\BracketBig{\psi^k} \AbsNormal{\nabla \widehat{\mu_\psi}}^2\right) \ \mathrm{d}x \\
    & \qquad + \frac12\BracketBig{\NormNormal{\nabla \phi}^2 + \beta \NormNormal{\nabla \psi}^2} + \int_\Omega\BracketNormal{F_{\phi}(\phi)+F_{\psi}(\psi)} \ \mathrm{d}x \\
    & \quad \le \frac{\lambda \alpha}{2} \int_\Omega \AbsBig{\sp^k}^2 \ \mathrm{d}x  + \frac12\BracketBig{\NormBig{\nabla \phi^k}^2+\beta \NormBig{\nabla \psi^k}^2} + \frac{\lambda \sigma_2}{2} \NormNormal{P_0 \phi^k}^2_* \\
    & \qquad + \int_\Omega\BracketBig{F_{\phi}\BracketBig{\phi^k}+F_{\psi}\BracketBig{\psi^k} +
   \lambda \AbsNormal{G(\phi,\psi)} + \lambda \AbsBig{G \BracketBig{\phi^k,\psi^k}}} \ \mathrm{d}x
    \nonumber \\
    & \qquad - \lambda h \BracketBig{\overline{\phi^k}-c} \int_\Omega\sigma_1\BracketBig{\phi^k} \widehat{\mu_\phi} \ \mathrm{d}x+ |\Omega| \BracketBig{ \overline{\phi} - \overline{\phi^k} } \BracketBig{ \overline{F'_\phi (\phi)} + \lambda \overline{G_\phi (\phi,\phi^k,\psi)} } \\
    & \quad \le C_k + 2 h |\Omega|^\frac{1}{2} \sigma_1^* \NormNormal{\widehat{\mu_\phi}} + |\Omega| \AbsBig{ \overline{\phi} - \overline{\phi^k} } \AbsBig{\overline{F'_\phi (\phi)} }.
\end{align*}
Recalling \cite[Section 3]{OGW2026}, the last two terms in the above inequality can be bounded by
\[
\frac14 h\underline{m_\phi} \|\nabla \widehat{\mu_\phi}\|^2 + \dfrac{4h \AbsNormal{\Omega} (\sigma_1^*)^2 C_P^2}{\underline{m_\phi}} + \frac14 h \underline{m_\phi} \|\nabla \widehat{\mu_\phi}\|^2 + \dfrac{ 4 |\Omega|^2 C_P^2 K_\phi^2 }{h\underline{m_\phi}} + 2 |\Omega| K_\phi.
\]
Therefore, it follows from the Poincar\'e--Wirtinger inequality that
\[
\NormNormal{\phi}_V + \NormNormal{\psi}_V + \NormNormal{\widehat{\mu_\phi}}_V + \NormNormal{\widehat{\mu_\psi}}_V \le C_k,
\]
where the positive constant $C_k$ depends on $k$ and $h$ but is uniform with respect to $\lambda \in [0,1]$.
So $\bm{w}$ belongs to a bounded set in $\widetilde{X}$ independent of $\lambda \in [0,1]$, and we obtain
\[
\NormNormal{\bm{f}}_{Y} = \NormNormal{\lambda \mathcal{F}_k (\bm{w})}_{Y} \le C_k,
\]
where the positive constant $C_k$ does not depend on $\lambda \in [0,1]$. By the Leray--Schauder principle, there exists a weak solution to the problem \eqref{dfdiscrete_project}--\eqref{chdiscrete4_project} that satisfies the discrete energy estimate inequality \eqref{discrete_energy_project} and possesses the additional regularity listed in Proposition \ref{discrete_exist_apriori}. The proof is complete.
\qed
\medskip

\begin{remark} \label{projection_to_origin}
\rm
As in \cite{OGW2026}, setting
    \[
    \mu_\phi = \widehat{\mu_\phi} + \overline{F'_\phi(\phi)} + \overline{G_\phi(\phi,\phi^k,\psi)},
    \quad
    \mu_\psi = \widehat{\mu_\psi} + \overline{F'_\psi(\psi)} + \overline{G_\psi(\phi^k,\psi,\psi^k)},
    \]
    we can easily verify that $(\sp,\phi,\psi,\mu_\phi,\mu_\psi)$ solves the problem
    \begin{numcases}{}
        \frac{\alpha}{h} (\sp-\sp^k,\bm{\theta}) + \InnerNormal{ \nu(\phi^k,\psi^k) \sp + \eta (\phi^k,\psi^k) |\sp|^{r-2} \sp }{ \bm{\theta} }\nonumber \\
        \quad = - \InnerNormal{ \phi^k \nabla \mu_\phi + \psi^k \nabla \mu_\psi }{ \bm{\theta} },
        & $\forall\, \bm{\theta} \in C^{\infty}_{0,\sigma} (\Omega)$,
        \label{dfdiscrete}\\
        \frac{\phi -\phi^k}{h} + \sp \cdot \nabla \phi^k + \sigma_1\BracketBig{\phi^k}
        \BracketBig{\overline{\phi^k} -c} = \mathrm{div} \BracketBig{m_\phi \BracketBig{\phi^k} \nabla \mu_\phi},
        & a.e. in $\Omega$,
        \label{chdiscrete1}\\
        \mu_\phi - G_{\phi} \BracketBig{\phi,\phi^k,\psi} - \sigma_2 \mathcal{N} \BracketBig{\phi - \overline{\phi}} = -\Delta \phi + F_{\phi}'(\phi),
        & a.e. in $\Omega$,
        \label{chdiscrete2}\\
        \frac{\psi -\psi^k}{h} + \sp \cdot \nabla \psi^k = \mathrm{div} \BracketBig{m_\psi \BracketBig{\psi^k} \nabla \mu_\psi},
        & a.e. in $\Omega$,
        \label{chdiscrete3}\\
        \mu_\psi - G_{\psi} \BracketBig{\phi^k,\psi,\psi^k} = -\beta \Delta \psi + F_{\psi}'(\psi),
        & a.e. in $\Omega$.
        \label{chdiscrete4}
    \end{numcases}
    Moreover, it satisfies the following energy inequality:
    \begin{align}
        & E_{\text{tot}} (\sp, \phi, \psi) + \frac{1}{2} \BracketBig{ \NormNormal{ \nabla \phi - \nabla \phi^k }^2 + \beta \NormNormal{ \nabla \psi - \nabla \psi^k }^2 } \nonumber \\
        & \qquad + \frac{\alpha}{2} \NormBig{ \sp - \sp^k }^2 + \frac{\sigma_2}{2} \NormBig{ \phi - \phi^k - \BracketBig{ \overline{\phi} - \overline{\phi^k} } }^2_* \nonumber \\
        & \qquad + h \int_\Omega \BracketBig{ \nu (\phi^k,\psi^k) | \sp |^2 + \eta (\phi^k,\psi^k) | \sp |^r + m_\phi (\phi^k) | \nabla \mu_\phi |^2 + m_\psi (\psi^k) | \nabla \mu_\psi |^2 } \ \mathrm{d} x \nonumber \\
        & \qquad + h \BracketBig{ \overline{\phi^k} - c } \int_\Omega \sigma_1(\phi^k) \mu_\phi \ \mathrm{d} x \nonumber \\
        & \quad \le E_{\text{tot}} (\sp^k, \phi^k, \psi^k). \label{discrete_energy}
    \end{align}
\end{remark}

\subsection{Proof of Theorem \ref{ws_chdf_thm}} \label{construct_weak_solution}

We now have all the ingredients to prove our first main result.
%
First, let us approximate the mobilities $m_\phi$ and $m_\psi$ by
\[
m_\phi^N \stackrel{\rm{def}}{=} \eta^N * m_\phi, \quad m_\psi^N \stackrel{\rm{def}}{=} \eta^N * m_\psi,
\]
where $\eta^N$ is a family of positive mollifiers in $\mathbb{R}$ and $*$ stands for the convolution in $\mathbb{R}$.
Then we have
\[
m_\phi^N \to m_\phi, \ m_\psi^N \to m_\psi \ \ \text{in } C(\mathbb{R}) \text{ as } N \to \infty.
\]
In addition, it holds
\[
m_\phi^N(s) \in \big[ \underline{m_\phi},m_\phi^* \big], \quad m_\psi^N(s) \in \big[ \underline{m_\psi},m_\psi^* \big],\quad
\forall\, s\in \mathbb{R},
\]
and
\[
\NormBig{ \BracketBig{m_\phi^N}' }_{C(\mathbb{R})} \le C_N\NormNormal{m_\phi}_{C(\mathbb{R})}, \quad  \NormBig{ \BracketBig{m_\psi^N}' }_{C(\mathbb{R})} \le C_N\NormNormal{m_\psi}_{C(\mathbb{R})},
\]
where the constant $C_N>0$ may depend on $N$, but we never use an explicit bound on $(m^N_\phi)'$ or on $(m^N_\psi)'$ in subsequent compactness arguments. This is only needed to fit the setting of Proposition  \ref{discrete_exist_apriori}.

Next, we approximate the initial data $(\phi_0,\psi_0) \in V \times V$ by $\phi_0^N,\psi_0^N \in W \times W$.
As shown in \cite{GP2022}, there exists a sequence $\big\{(\phi_0^N,\psi_0^N)\big\}_{N=1}^{\infty} \subset H^3_N(\Omega) \times H^3_N(\Omega)$ satisfying
\begin{equation*}
    \left\|\phi_0^N \right\|_{L^{\infty}(\Omega)},\ \left\|\psi_0^N \right\|_{L^{\infty}(\Omega)} \le 1- \frac{1}{N}, \quad  \forall\, N \in \mathbb{N_+},
\end{equation*}
and
\begin{equation*}
    \BracketBig{\phi_0^N,\psi_0^N} \to (\phi_0,\psi_0) \ \ \mathrm{in} \ \ V \ \ \mathrm{as} \ \ N \to \infty,
\end{equation*}
where $H^3_N(\Omega)\stackrel{\rm{def}}{=} W \cap H^3(\Omega)$.
Moreover, an argument similar to that in \cite[Remark 3.2]{OGW2026} yields, if $N \ge 1+ \sigma_1^*$, then $\phi^{k+1}$ satisfies
\begin{align*}
    \AbsMax{ \overline{\phi^{k+1}} - c } & =  \prod\limits_{n=0}^k \BracketBig{1-\frac{\overline{\sigma_1\BracketNormal{\phi^n}}}{N}} \AbsMax{\overline{\phi^N_0} - c} \\
    & \le e^{- \frac{1}{N} \sum\limits_{n=0}^k \overline{\sigma_1\BracketNormal{\phi^n}}} \AbsMax{\overline{\phi^N_0} - c} \le 2 e^{- \frac{1}{N} \sum\limits_{n=0}^k \overline{\sigma_1\BracketNormal{\phi^n}}}.
\end{align*}

Now we can iteratively apply Proposition \ref{discrete_exist_apriori} combined with Remark \ref{projection_to_origin} to construct a sequence of solutions $\SetBig{ \sp^{k+1,N}, \phi^{k+1,N}, \psi^{k+1,N}, \mu_\phi^{k+1,N}, \mu_\psi^{k+1,N} }_{k=0}^{\infty}$.
In the subsequent analysis, the following convention will be used. For a given sequence $\SetNormal{f^{k,N}}_{k=0}^{\infty}$ that may depend on $N$, we denote by $f^N$ the function defined in $[-h,\infty)$ such that
\[
f^N(t) = f^{k,N}, \ \ \  t \in [(k-1)h,kh),\quad \forall\, k \in \mathbb{N}.
\]
In addition, we use $\widetilde{f}^N$ to denote the piecewise linear interpolant of $f^N(t^{k,N})$, where $t^{k,N} = kh, \ k \in \mathbb{N}$, that is,
\[
\widetilde{f}^N = \frac{1}{h} \chi_{[0,h]} *_t f^N.
\]
For a given function $f=f(t)$, we set
\begin{align*}
    & (\Delta^+_h f)(t) \stackrel{\rm{def}}{=} f(t+h)-f(t), \qquad (\Delta^-_h f)(t) \stackrel{\rm{def}}{=} f(t)-f(t-h),\\
    & \partial^+_{t,h} f(t) \stackrel{\rm{def}}{=} \frac{1}{h} (\Delta^+_h f)(t), \qquad\qquad \ \,\partial^-_{t,h} f(t) \stackrel{\rm{def}}{=} \frac{1}{h} (\Delta^-_h f)(t),\\
    & f_h \stackrel{\rm{def}}{=} (\tau^*_h f)(t) = f(t-h).
\end{align*}
For any $\bm{\theta} \in C^\infty_0 (Q)^d$ satisfying $\mathrm{div} \ \bm{\theta} =0$, we choose
\[
\bm{\theta}^k = \int_{kh}^{(k+1)h} \theta(t) \ \mathrm{d} t
\]
as the test function in \eqref{dfdiscrete} and sum over $k\in \mathbb{N}$. This gives
\begin{align}
    & \int_0^\infty \InnerBig{ \partial^-_{t,h} \sp^N }{\bm{\theta}} \ \mathrm{d}t + \int_0^\infty \InnerBig{ \nu (\phi^N_h,\psi^N_h) \sp^N }{\bm{\theta}} \ \mathrm{d}t + \int_0^\infty \PairingBig{ \eta (\phi^N_h,\psi^N_h) \AbsBig{\sp^N}^{r-2} \sp^N }{\bm{\theta}}_{\bm{L}^{r'}(\Omega),\bm{L}^r(\Omega)} \ \mathrm{d}t \nonumber \\
    & \quad = - \int_0^\infty \InnerBig{ \phi^N_h \nabla \mu_\phi^N + \psi^N_h \nabla \mu_\psi^N }{\bm{\theta}} \ \mathrm{d} t. \label{dfglobaldiscrete}
\end{align}
Analogously, we obtain
\begin{align}
    & - \int_0^\infty \InnerBig{\partial^-_{t,h} \phi^N}{\zeta} \ \mathrm{d} t + \int_0^\infty  \InnerBig{\sp^N \phi^N_h}{\nabla \zeta} \ \mathrm{d}t - \int_0^\infty \BracketBig{ \overline{\phi^N_h} - c } \InnerBig{\sigma_1\BracketBig{\phi^N_h}}{\zeta} \ \mathrm{d} t \nonumber \\
    & \quad = \int_0^\infty \InnerBig{m_\phi^N \BracketBig{\phi^N_h} \nabla \mu_\phi^N}{\nabla \zeta}\ \mathrm{d} t,\label{chglobaldiscrete1}\\
    & -\int_0^\infty \InnerBig{\partial^-_{t,h} \psi^N}{\zeta}\ \mathrm{d} t + \int_0^\infty \InnerBig{\sp^N \psi^N_h}{\nabla \zeta}\ \mathrm{d} t= \int_0^\infty \InnerBig{m_\psi^N \BracketBig{\psi^N_h} \nabla \mu_\psi^N}{\nabla \zeta}\ \mathrm{d} t,\label{chglobaldiscrete3}
\end{align}
for all $\zeta \in C^{\infty}_0\big((0,\infty);C^1\big(\overline{\Omega} \big)\big)$, and
\begin{align}
    & \mu_\phi^N - G_{\phi} \BracketBig{\phi^N,\phi^N_h,\psi^N}
    -\sigma_2 \mathcal{N} \BracketBig{ \phi^N - \overline{\phi^N} }
    = -\Delta \phi^N + F_{\phi}' \BracketBig{\phi^N}, \label{chglobaldiscrete2}\\
    & \mu_\psi^N - G_{\psi} \BracketBig{\phi^N_h,\psi^N,\psi^N_h} = -\beta \Delta \psi^N + F_{\psi}' \BracketBig{\psi^N}, \label{chglobaldiscrete4}
\end{align}
almost everywhere in $Q$.
On the other hand, it follows from the discrete energy inequality \eqref{discrete_energy} that for $\widetilde{E}^N$, i.e. the linear interpolant of $E^N(t)=E_{\text{tot}} \BracketBig{\sp^{k,N},\phi^{k,N},\psi^{k,N}}$, $k \in \mathbb{N}$, we have
\begin{align*}
    -\frac{\mathrm{d}}{\mathrm{d}t} \widetilde{E}^N (t) & \ge \int_\Omega \BracketBig{ \nu (\phi^N_h(t),\psi^N_h(t)) \AbsBig{\sp^N(t)}^2 + \eta (\phi^N_h(t),\psi^N_h(t)) \AbsBig{\sp^N(t)}^r }  \ \mathrm{d} x \\
    & \quad + \int_\Omega \BracketBig{ m_\phi^N \BracketBig{\phi^N_h(t)} \AbsBig{\nabla \mu_\phi^N(t)}^2 + m_\psi^N \BracketBig{\psi^N_h(t)} \AbsBig{\nabla \mu_\psi^N(t)}^2 }  \ \mathrm{d} x \\
    & \quad + \int_\Omega \BracketBig{ \overline{\phi^N_h}(t) -c } \sigma_1 \BracketBig{\phi^N_h(t)} \mu_\phi^N(t) \ \mathrm{d} x,
\end{align*}
for all $t \in (t^{k,N},t^{k+1,N})$, $k \in \mathbb{N}$, which implies
\begin{align}
    & \widetilde{E}^N (t) + \int_{Q_t} \BracketBig{ \nu (\phi^N_h,\psi^N_h) \AbsBig{\sp^N}^2 + \eta (\phi^N_h,\psi^N_h) \AbsBig{\sp^N}^r }  \ \mathrm{d} (x,\tau) \nonumber \\
    & \qquad + \int_{Q_t} \BracketBig{ m_\phi^N \BracketBig{\phi^N_h} \AbsBig{\nabla \mu_\phi^N}^2 + m_\psi^N \BracketBig{\psi^N_h} \AbsBig{\nabla \mu_\psi^N}^2 }  \ \mathrm{d} (x,\tau) \nonumber \\
    & \quad \le \widetilde{E}^N(0) - \int_{Q_t} \BracketBig{ \overline{\phi^N_h} -c } \sigma_1 \BracketBig{\phi^N_h} \mu_\phi^N \ \mathrm{d} (x,\tau),
    \quad \forall\, t\in [0,\infty). \label{energy_global_approximate}
\end{align}
By choosing a positive constant $\delta\in (0,\frac12)$ such that $c, \overline{\phi_0} \in [-1+\delta,1-\delta]$ and $\overline{\psi_0} \in [\delta,1-\delta]$, we have $\overline{\phi^N} \in [-1+\delta,1-\delta]$ and $\overline{\psi^N} \in [\delta,1-\delta]$.
Thus, using Lemma \ref{apriori_ch}, we can conclude that there exist two positive constants $K_\phi = K_\phi(\delta)$ and $K_\psi=K_\psi(\delta)$ such that
\[
\NormBig{ \mu_\phi^N } \le K_\phi(\delta) \BracketBig{ 1 + \NormBig{\nabla \mu_\phi^N} }, \qquad  \NormBig{ \mu_\psi^N } \le K_\psi(\delta) \BracketBig{ 1 + \NormBig{\nabla \mu_\psi^N} },
\]
for all $N \ge 1 + \sigma_1^*$.
Therefore, the following estimate holds
\begin{equation*}
    \AbsMax{ \int_{\Omega} \BracketBig{ \overline{\phi^N_h} -c } \sigma_1 \BracketBig{\phi^N_h} \mu_\phi^N \ \mathrm{d} x }
    \le \frac{1}{4} \underline{m_\phi} \NormNormal{\nabla \mu_\phi^N}^2 + \dfrac{4 |\Omega| (\sigma_1^*)^2 K_\phi (\delta)^2}{\underline{m_\phi}} + 2 K_\phi (\delta) \sigma_1^* \AbsNormal{\Omega}^{\frac{1}{2}}.
\end{equation*}
Then, for any given $T>0$, we deduce from \eqref{energy_global_approximate} that the following terms are uniformly bounded with respect to $N$, in the corresponding spaces (cf. \cite{OGW2026,ADG2013ndm}):
\begin{align*}
    & \sp^N \in L^r \BracketBig{ 0,T; \bm{L}^r_\sigma(\Omega) }, \quad  \sqrt\alpha \sp^N \in L^\infty \BracketBig{0,T;\bm{L}^2_\sigma(\Omega)}, \\
    & \phi^N, \psi^N \in L^\infty \BracketBig{ 0,T;V },\quad \mu_\phi^N, \mu_\psi^N \in L^2 \BracketBig{ 0,T; V }.
\end{align*}
Using a diagonal argument as in \cite{BF2013} for the classical Navier--Stokes equations for incompressible viscous fluids, we can extract a convergent subsequence and find a limit $(\sp,\phi,\psi,\mu_\phi,\mu_\psi)$, defined in $[0,\infty)$, which satisfies
\begin{align}
    & \sp^N \rightharpoonup \sp\ \text{in}\ L^r \BracketBig{ 0,T; \bm{L}^r_\sigma(\Omega) },\ \  \sqrt \alpha \sp^N \stackrel{*}{\rightharpoonup} \sqrt \alpha \sp\ \text{in}\  L^\infty \BracketBig{0,T;\bm{L}^2_\sigma(\Omega)},
    \notag\\
    &\phi^N \stackrel{*}{\rightharpoonup} \phi\ \ \text{and}\ \ \psi^N \stackrel{*}{\rightharpoonup} \psi\ \text{in}\ L^\infty \BracketBig{ 0,T; V },
    \notag\\
    & \mu_\phi^N \rightharpoonup \mu_\phi\ \ \text{and}\ \ \mu_\psi^N \rightharpoonup \mu_\psi \ \text{in}\ L^2 \BracketBig{ 0,T; V },
    \notag
\end{align}
for any $T>0$. Hereafter, all limits are meant to be for a suitable subsequence $N_k \to \infty$ (resp. $h_k \to 0$) for $k \to \infty$ unless otherwise stated.

Arguing as in \cite[Section 3.2]{OGW2026} (see also \cite{ADG2013ndm}), for $\BracketBig{\widetilde{\phi}^N,\widetilde{\psi}^N}$, we get
\[
\widetilde{\phi}^N \to \phi, \ \widetilde{\psi}^N \to \psi \ \mathrm{in} \ C\BracketBig{ [0,T];H }, \quad \forall\, T>0.
\]
Furthermore, we can check that
\begin{align*}
    & \phi^N \to \phi, \ \psi^N \to \psi \ \mathrm{in} \ L^2 \BracketBig{ 0,T; V }, \\
    & \phi, \psi \in C_w \BracketBig{[0,\infty);V} \cap H^1 \BracketBig{0,T;V^*} \cap L^2 \BracketBig{0,T;W^{2,p}(\Omega)}, \\
    & (\phi(0), \psi(0)) = (\phi_0,\psi_0), \\
    & \partial \widetilde{E}_\phi (\phi) = \mu_\phi - \partial_\phi G(\phi,\psi) - \sigma_2 \mathcal{N} \BracketBig{\phi - \overline{\phi}}, \\
    & \partial \widetilde{E}_\psi (\psi) = \mu_\psi - \partial_\psi G(\phi,\psi),
\end{align*}
with $p=6$ if $d=3$ and for any $p \in [2,\infty)$ if $d=2$.
These convergence results enable us to pass to the limit in equations \eqref{chglobaldiscrete1}--\eqref{chglobaldiscrete4}.

Next, we consider the limit for equation \eqref{dfglobaldiscrete}.
Firstly, because of the uniform boundedness of $\eta (\phi^N_h,\psi^N_h)^{\tfrac{r-1}{r}} \AbsBig{\sp^N}^{r-2} \sp^N $ in $L^{r'} \BracketBig{0,T;\bm{L}^{r'}(\Omega)}$ with respect to $N$ for any $T>0$, it follows from a diagonal argument that there exists some $\Psi \in L^{r'}_{\mathrm{loc}} \BracketBig{[0,\infty);\bm{L}^{r'}(\Omega)}$ such that, up to a subsequence,
\[
\eta (\phi^N_h,\psi^N_h)^\frac{r-1}{r} \AbsBig{\sp^N}^{r-2} \sp^N \rightharpoonup \Psi \ \ \mathrm{in} \ L^{r'} \BracketBig{0,t;\bm{L}^{r'}(\Omega)}, \quad \forall\, t>0.
\]
 Due to the strong convergence of $\phi^N$ and $\psi^N$, we also find
\[
\eta (\phi^N_h,\psi^N_h) \AbsBig{\sp^N}^{r-2} \sp^N \rightharpoonup \eta (\phi,\psi)^\frac{1}{r} \Psi \ \ \mathrm{in} \ L^{r'} \BracketBig{0,t;\bm{L}^{r'}(\Omega)}.
\]
Then, the limit $(\sp,\phi,\psi,\mu_\phi,\mu_\psi)$ solves the following problem:
\begin{numcases}{}
    - \alpha \InnerNormal{\sp}{\partial_t \bm{v}} + \PairingNormal{ \eta (\phi,\psi)^\frac{1}{r} \Psi }{ \bm{v} }_{ L^{r'} \BracketNormal{ 0,T;\bm{L}^{r'}(\Omega) },L^r \BracketNormal{ 0,T;\bm{L}^r(\Omega) } }  \nonumber \\
    \quad + \InnerNormal{\nu (\phi,\psi) \sp}{\bm{v}} = - \InnerNormal{ \phi \nabla \mu_\phi + \psi \nabla \mu_\psi }{ \bm{v} } \label{eq:df1*}, & $\forall\,\bm{v} \in C^\infty_0 \BracketBig{(0,T);\bm{L}^r_\sigma(\Omega)}$,\\
    \PairingNormal{ \partial_t \phi + \sigma_1(\phi) \BracketBig{\overline{\phi} -c} }{ \zeta }_{L^2 (0,T;V^*),L^2(0,T;V)} - \InnerNormal{ \phi \sp }{\nabla \zeta}_{Q_T} \nonumber \\
    \quad = - \InnerNormal{ m_{\phi} (\phi) \nabla \mu_\phi }{ \nabla \zeta }_{Q_T} \label{eq:ch1*}, & $\forall\, \zeta \in L^2\BracketBig{0,T;V}$, \\
    \mu_{\phi} = -\Delta \phi + \sigma_2 \mathcal{N} \BracketBig{\phi - \overline{\phi}} + F_{\phi}'(\phi) + \partial_{\phi} G(\phi,\psi) \label{eq:ch2*}, & a.e. in $Q_T$,\\
    \PairingNormal{ \partial_t \psi }{ \zeta }_{L^2 (0,T;V^*),L^2(0,T;V)} - \InnerNormal{ \psi \sp }{\nabla \zeta}_{Q_T} = - \InnerNormal{ m_{\psi} (\psi) \nabla \mu_\psi }{ \nabla \zeta }_{Q_T} \label{eq:ch3*}, & $\forall\, \zeta \in L^2\BracketBig{0,T;V}$, \\
    \mu_{\psi} = - \beta \Delta \psi + F_{\psi}'(\psi) + \partial_{\psi} G(\phi,\psi) \label{eq:ch4*}, & a.e. in $Q_T$,
\end{numcases}
for any given $T>0$, subject to the following initial conditions:
\begin{numcases}{}
    \sp|_{t=0} = \sp_0, & in $\Omega$, if $\alpha >0$, \nonumber \\
    \phi|_{t=0} = \phi_0(x), \ \psi|_{t=0} = \psi_0(x),  & in $\Omega$.\label{eq:chi*}
\end{numcases}

It remains to verify that $\Psi = \eta (\phi,\psi)^\frac{r-1}{r} \AbsBig{\sp}^{r-2} \sp$.
We note that the regularity property $\phi, \psi \in H^1(0,T; V^*)$ combined with a standard argument (see, e.g., \cite[Lemma 4.1]{RS2004}) ensures the absolute continuity of the energy components $E_\phi(\phi)$ and $E_\psi(\psi)$ on $[0,T]$.
In addition, due to the weak lower semicontinuity of the convex potential functions, the continuity of the energy components ensures that $\NormNormal{\phi(t)}_V$ and $\NormNormal{\psi(t)}_V$ are continuous in time.
Combining the norm continuity with the weak continuity of $\phi$ and $\psi$ with values in $V$, we readily deduce the strong continuity (see \cite[Section 3]{GGW2018} and \cite[Theorem 6]{A2009}):
\[
\phi,\psi \in BC\BracketBig{[0,\infty);V}.
\]
Furthermore, the classical result on the Gelfand triple (see \cite[Proposition 23.23]{Z1990}) implies the absolute continuity of $\NormNormal{\sp}^2$ in time whenever $\partial_t \sp \in L^{r'}_{\mathrm{loc}}(0,\infty;\bm{L}^{r'}_\sigma(\Omega))$.
These absolute continuity properties are crucial for rigorously establishing the energy identity.

In the following, we distinguish two cases.
If $\alpha = 0$, using a standard argument (see \cite[Lemma 4.1]{RS2004}), we can take $\bm{v}=\bm{u}$ in \eqref{eq:df1*}, $\zeta= \mu_\phi$ in \eqref{eq:ch1*} and $\zeta = \mu_\psi$ in \eqref{eq:ch3*} to get
\begin{align*}
    & E_{\text{free}} (\phi(t),\psi(t)) + \int_{Q_t} \BracketBig{ m_\phi (\phi) \AbsNormal{ \nabla \mu_\phi }^2 + m_\psi (\psi) \AbsNormal{ \nabla \mu_\psi }^2 + \nu (\phi,\psi) \AbsNormal{\sp}^2 } \ \mathrm{d} (x,\tau) \\
    & \qquad + \PairingNormal{ \Psi }{ \eta (\phi,\psi)^\frac{1}{r} \sp }_{L^{r'}(0,t;\bm{L}^{r'}(\Omega)),L^r(0,t;\bm{L}^r(\Omega))} \\
    & \quad = E_{\text{free}} (\phi_0,\psi_0) - \int_{Q_t} (\overline{\phi} - c) \sigma_1(\phi) \mu_\phi \ \mathrm{d} (x,\tau),\quad \forall\, t\geq 0.
\end{align*}
On the other hand, the strong convergence of $\phi^N$ and $\psi^N$ ensures that
\[ \liminf\limits_{N \to \infty} \widetilde{E}^N (t) \ge E_{\text{free}} ( \phi(t),\psi(t) ),\quad \text{for a.a}\ t\in [0,\infty),
\]
and
\begin{align*}
& \sqrt{\nu (\phi^N_h,\psi^N_h)} \sp^N \rightharpoonup \sqrt{\nu (\phi,\psi)} \sp, \ \mathrm{in} \ L^2\BracketBig{0,T;L^2(\Omega)}, \\
&\eta (\phi^N_h,\psi^N_h)^\frac{1}{r} \sp^N \rightharpoonup \eta (\phi,\psi)^\frac{1}{r} \sp \ \mathrm{in} \ L^r\BracketBig{0,T;\bm{L}^r(\Omega)} ,
\end{align*}
for any given $T>0$.
Hence, we have
\begin{align*}
    & \limsup\limits_{N \to \infty} \int_{Q_t}  \eta (\phi^N_h,\psi^N_h) \AbsBig{\sp^N}^r  \ \mathrm{d} (x,\tau) \\
    & \quad \le \limsup\limits_{N \to \infty} E_{\text{free}} (\phi^N_0,\psi^N_0) - \liminf\limits_{N \to \infty} \widetilde{E}^N(t) - \liminf\limits_{N \to \infty} \int_{Q_t} \nu (\phi^N_h,\psi^N_h) \AbsBig{\sp^N}^2  \ \mathrm{d} (x,\tau) \\
    & \qquad - \liminf\limits_{N \to \infty} \int_{Q_t} \BracketBig{ m_\phi^N (\phi^N_h) \AbsBig{\nabla \mu_\phi^N}^2 + m_\psi^N (\psi^N_h) \AbsBig{\nabla \mu_\psi^N}^2 } \ \mathrm{d} (x,\tau) \\
    & \qquad - \liminf\limits_{N \to \infty} \int_{Q_t} (\overline{\phi^N_h}-c) \sigma_1 (\phi^N_h) \mu_\phi^N \ \mathrm{d} (x,\tau) \\
    & \quad \le E_{\text{free}} (\phi_0,\psi_0) - E_{\text{free}}(\phi(t),\psi(t)) - \int_{Q_t} \nu (\phi,\psi) \AbsNormal{\sp}^2 \ \mathrm{d} (x,\tau) \\
    & \qquad - \int_{Q_t} \BracketBig{ m_\phi (\phi) \AbsBig{\nabla \mu_\phi}^2 + m_\psi (\psi) \AbsBig{\nabla \mu_\psi}^2 } \ \mathrm{d} (x,\tau)
    - \int_{Q_t} (\overline{\phi}-c) \sigma_1 (\phi) \mu_\phi \ \mathrm{d} (x,\tau) \\
    & \quad = \PairingNormal{ \Psi }{ \eta (\phi,\psi)^\frac{1}{r} \sp }_{L^{r'}(0,t;\bm{L}^{r'}(\Omega)),L^r(0,t;\bm{L}^r(\Omega))},
\end{align*}
for all $t \in \SetNormal{\tau \ge 0: \lim\limits_{N \to \infty} \widetilde{E}^N (\tau) = E_{\text{free}}(\phi(\tau),\psi(\tau))}$.
Here we used the lower semicontinuity properties of norms thanks to the weak convergence of $\mu_\phi^N$ and $\mu_\psi^N$.
Therefore, by the monotonicity and hemicontinuity of the operator $\mathcal{M}_t:\bm{L}^r(Q_t) \to \bm{L}^{r'}(Q_t)$, with
\[
\mathcal{M}(\bm{f}) = \AbsNormal{\bm{f}}^{r-2} \bm{f}, \quad  \forall\, \bm{f} \in \bm{L}^r(Q_t),
\]
we infer from the abstract result \cite[Theorem 9.13-2(a)]{Ciarlet2013} (see also the references therein) that
\[ \Psi = \mathcal{M} (\eta (\phi,\psi)^\frac{1}{r} \sp) \ \ \mathrm{in} \ \ L^{r'}(Q_t), \quad \text{for a.a.}\ t\in[0,\infty).
\]
Hence, we see that the limit $\BracketNormal{ \sp,\phi,\psi,\mu_\phi,\mu_\psi }$ is a solution to the problem \eqref{eq:chdf1}--\eqref{eq:chdfi}.

Let us consider now the case $\alpha >0$. For any $T>0$, we have
\[
\partial_t \widetilde{\sp}^N \in L^{r'} \BracketBig{0,T;\bm{L}^{r'}_\sigma(\Omega)},
\]
is uniformly bounded with respect to $N$.
Then, a standard argument, combined with the regularity of $\sp$ inherited from $\widetilde{\sp}^N$ (namely, $\sp \in W^{1,r'}_{\mathrm{loc}}\BracketBig{0,\infty;\bm{L}^{r'}_\sigma(\Omega)} \cap L^r_{\mathrm{loc}}\BracketBig{0,\infty;\bm{L}^r_\sigma(\Omega)}$), yields
\begin{align}
    & \PairingNormal{ \alpha \partial_t \sp + \eta (\phi,\psi)^\frac{1}{r} \Psi }{ \bm{v} }_{ L^{r'}(0,T;\bm{L}^{r'}(\Omega)),L^r(0,T;\bm{L}^r(\Omega)) }  \nonumber \\
    & \quad + \InnerNormal{\nu (\phi,\psi) \sp}{\bm{v}} = - \InnerNormal{ \phi \nabla \mu_\phi + \psi \nabla \mu_\psi }{ \bm{v} }, \quad \forall\, \bm{v} \in L^r \BracketBig{ 0,T;\bm{L}^r_\sigma(\Omega) }. \label{eq:df1**}
\end{align}
Like in the case $\alpha = 0$, applying \cite[Lemma 4.1]{RS2004} and \cite[Proposition 23.23]{Z1990}, we can obtain
the energy identity
\begin{align*}
    & E_{\text{tot}} (\sp(t),\phi(t),\psi(t)) + \int_{Q_t} \BracketBig{ m_\phi (\phi) \AbsNormal{ \nabla \mu_\phi }^2 + m_\psi (\psi) \AbsNormal{ \nabla \mu_\psi }^2 + \nu (\phi,\psi) \AbsNormal{\sp}^2 } \ \mathrm{d} (x,\tau) \\
    & \qquad + \PairingNormal{ \Psi }{ \eta (\phi,\psi)^\frac{1}{r} \sp }_{L^{r'}(0,t;\bm{L}^{r'}(\Omega)),L^r(0,t;\bm{L}^r(\Omega))} \\
    & \quad = E_{\text{tot}} (\sp_0,\phi_0,\psi_0) - \int_{Q_t} (\overline{\phi} - c) \sigma_1(\phi) \mu_\phi \ \mathrm{d} (x,\tau),\quad \forall\, t\geq 0.
\end{align*}
In addition, from equation \eqref{dfglobaldiscrete} it follows that
\[ \NormBig{ \widetilde{\sp}^N - \sp^N }_{\bm{L}^{r'}_\sigma(\Omega)} \le h \NormBig{ \partial_t \widetilde{\sp}^N }_{\bm{L}^{r'}_\sigma(\Omega)}, \]
which implies
\begin{align}
&\NormBig{ \widetilde{\sp}^N(t) - \sp^N(t) }_{\bm{L}^{r'}_\sigma(\Omega)} \le h^{\frac{1}{r}} \NormBig{\partial_t \widetilde{\sp}^N}_{L^{r'} \BracketBig{0,T;\bm{L}^{r'}_\sigma(\Omega)}},
\quad \forall\, t\in [0,T],
\label{ww}\\
& \widetilde{\sp}^N \rightharpoonup \sp \ \mathrm{in} \ W^{1,r'} \BracketBig{ 0,T;\bm{L}^{r'}_\sigma(\Omega) } \cap L^r \BracketBig{0,T;\bm{L}^r_\sigma(\Omega)} \hookrightarrow C \BracketBig{[0,T];\bm{L}^2_\sigma (\Omega)}.
\label{weakx}
\end{align}
Since $T>0$ is arbitrary, we can deduce from \eqref{weakx} that
\begin{align}
 \sp^N(t) \rightharpoonup \sp(t) \ \mathrm{in} \  \bm{L}^2_\sigma (\Omega),
 \label{weakt}
 \end{align}
for all $t \in [0,\infty)$ including $t = 0$. To conclude \eqref{weakt}, we proceed as follows.
Using the regularity of $\sp,\widetilde{\sp}^N$, we can apply the abstract result \cite[Lemma 7.3]{Rubi} to get the identity
\[
\langle \widetilde{\sp}^N(t)-{\sp}(t), {\bm v}\rangle
= \int_0^t \langle \partial_t (\widetilde{\sp}^N -{\sp})(s), {\bm v}\rangle\, \mathrm{d}s = \int_0^T \langle \partial_t (\widetilde{\sp}^N -{\sp})(s), {\bm v}\chi_{[0,t]}\rangle\,\mathrm{d}s,
\]
for any fixed $t\in [0,T]$ and any ${\bm v} \in \bm{L}^{r}(\Omega)$. Thus, from the weak convergence to $\mathbf 0$ of $\{\partial_t (\widetilde{\sp}^N -{\sp})\}$ in $L^{r'} \BracketBig{ 0,T;\bm{L}^{r'}_\sigma(\Omega) } $ (see \eqref{weakx}),
we infer that
\[
\widetilde{\sp}^N(t) \rightharpoonup \sp(t) \quad \text{weakly in } \bm{L}^{r'}_\sigma(\Omega), \quad \forall\, t\in [0,T].
\]
Since $\widetilde \sp_N$ is uniformly bounded in $C([0,T];\bm{L}^2_\sigma(\Omega))$, a density argument gives the convergence of $\widetilde \sp _N$, which, combined with \eqref{ww}, leads to \eqref{weakt}.

The above convergence results imply that
\[ \liminf\limits_{N \to \infty} \widetilde{E}^N (t) \ge E_{\text{tot}} (\sp(t),\phi(t),\psi(t)),
\quad \text{for a.a.}\ t \in [0,\infty).
\]
Therefore, it holds
\begin{align*}
    & \limsup\limits_{N \to \infty} \int_{Q_t}  \eta (\phi^N_h,\psi^N_h) \AbsBig{\sp^N}^r  \ \mathrm{d} (x,\tau) \\
    & \quad \le \limsup\limits_{N \to \infty} E_{\text{tot}} (\sp_0,\phi^N_0,\psi^N_0) - \liminf\limits_{N \to \infty} \widetilde{E}^N(t) - \liminf\limits_{N \to \infty} \int_{Q_t} \nu (\phi^N_h,\psi^N_h) \AbsBig{\sp^N}^2 \ \mathrm{d} (x,\tau) \\
    & \qquad - \liminf\limits_{N \to \infty} \int_{Q_t} \BracketBig{ m_\phi^N (\phi^N_h) \AbsBig{\nabla \mu_\phi^N}^2 + m_\psi^N (\psi^N_h) \AbsBig{\nabla \mu_\psi^N}^2 } \ \mathrm{d} (x,\tau) \\
    & \qquad - \liminf\limits_{N \to \infty} \int_{Q_t} (\overline{\phi^N_h}-c) \sigma_1 (\phi^N_h) \mu_\phi^N \ \mathrm{d} (x,\tau) \\
    & \quad \le E_{\text{tot}} (\sp_0,\phi_0,\psi_0) - E_{\text{tot}}(\sp(t),\phi(t),\psi(t)) - \int_{Q_t} \nu (\phi,\psi) \AbsNormal{\sp}^2 \ \mathrm{d} (x,\tau) \\
    & \qquad - \int_{Q_t} \BracketBig{ m_\phi (\phi) \AbsBig{\nabla \mu_\phi}^2 + m_\psi (\psi) \AbsBig{\nabla \mu_\psi}^2 } \ \mathrm{d} (x,\tau)
    - \int_{Q_t} (\overline{\phi}-c) \sigma_1 (\phi) \mu_\phi \ \mathrm{d} (x,\tau) \\
    & \quad = \PairingNormal{ \Psi }{ \eta (\phi,\psi)^\frac{1}{r} \sp }_{L^{r'}(0,t;\bm{L}^{r'}(\Omega)),L^r(0,t;\bm{L}^r(\Omega))},
\end{align*}
for all $t \in \SetNormal{\tau \ge 0: \liminf\limits_{N \to \infty} \widetilde{E}^N (\tau) \ge E_{\text{tot}}(\sp(\tau),\phi(\tau),\psi(\tau)) }$.
Arguing as for the case $\alpha=0$, we find
\[ \Psi = \eta (\phi,\psi)^\frac{r-1}{r} \AbsNormal{\sp}^{r-2} \sp. \]
Thus, $(\sp,\phi,\psi,\mu_\phi,\mu_\psi)$ is a weak solution to problem \eqref{eq:chdf1}--\eqref{eq:chdfi} in the sense of Definition \ref{ws_chdf_def}.

To recover the pressure $\pi$, we define an auxiliary function $\bm{f}_{\alpha}$ belonging to $W^{1,r'}_{\mathrm{loc}} \BracketBig{ [0,\infty);\bm{L}^{r'}(\Omega) } \subset C \BracketBig{ [0,\infty);\bm{L}^{r'}(\Omega) }$ by setting
\[
\bm{f}_{\alpha}(t) \stackrel{\rm{def}}{=} \int_0^t \big(\mu_\phi \nabla \phi + \mu_\psi \nabla \psi - \alpha \partial_t \sp - \nu (\phi,\psi) \sp - \eta (\phi,\psi) \AbsNormal{\sp}^{r-2} \sp \big)\, \mathrm{d} \tau.
\]
Here, we used the property
\[\phi,\psi \in L^\infty(Q) \cap L^4_{\mathrm{loc}} \BracketBig{ 0,\infty;W } \cap L^2_{\mathrm{loc}} \BracketBig{ 0,\infty;W^{2,6}(\Omega) }\]
combined with the Gagliardo--Nirenberg--Sobolev inequality in the three-dimensional case, while for the two-dimensional case we use $ \phi,\psi \in L^\infty(0,T;V) $ for all $T>0$. Indeed, we have
$\mu_\phi \nabla \phi + \mu_\psi \nabla \psi \in L^{\frac{2}{2-\theta}}_{\mathrm{loc}} \BracketBig{ [0,\infty);L^{\frac{2}{\theta}}(\Omega) }$ for all $\theta \in (0,1)$ if $d=3$ and $\mu_\phi \nabla \phi + \mu_\psi \nabla \psi \in L^2_{\mathrm{loc}} \BracketBig{ [0,\infty);L^q(\Omega) }$ for any $q \in (1,2)$ if $d=2$.
Since $r' \in (1,2)$, we can conclude that $\mu_\phi \nabla \phi + \mu_\psi \nabla \psi \in L^{r'}_{\mathrm{loc}} \BracketBig{ [0,\infty);\bm{L}^{r'}(\Omega) }$.
Applying \cite[Theorem 4.2]{BS1990}, we see that for all $t \in [0,\infty)$, there exists $p(t) \in L^{r'}_{(0)}(\Omega)$ satisfying
\[
\nabla p(t) = \bm{f}_{\alpha}(t) \in \bm{L}^{r'}(\Omega).
\]
This gives $p \in C\BracketBig{[0,\infty);W^{1,r'}_{(0)} (\Omega)}$. In addition, it holds
\begin{align*}
    \int_0^T \PairingNormal{p}{\partial_t g}_{W^{1,r'}_{(0)}(\Omega), (W^{1,r'}_{(0)}(\Omega))^*} \ \mathrm{d} t
     &= - \int_0^T \PairingNormal{ \partial_t \bm{f}_{\alpha} }{ \nabla \mathcal{N} g }_{\bm{L}^{r'}(\Omega),\bm{L}^r(\Omega)} \ \mathrm{d} t\\
    &\le C \NormNormal{\bm{f}_{\alpha}}_{W^{1,r'}(0,T;\bm{L}^{r'}(\Omega))} \NormNormal{g}_{L^r( 0,T;(W^{1,r'}_{(0)}(\Omega))^* )} ,
\end{align*}
for all $g \in C^\infty_0 \BracketBig{ (0,T);(W^{1,r'}_{(0)}(\Omega))^* }$.
Thus, there exists $\pi \in L^{r'}_{\mathrm{loc}} \BracketBig{ 0,\infty;W^{1,r'}_{(0)} (\Omega) } $ such that $\pi = \partial_t p $.
In particular, we have $\nabla \pi = \partial_t \bm{f}_{\alpha}$, that is,
\[ \alpha \partial_t \sp + \nu (\phi,\psi) \sp + \eta (\phi,\psi) \AbsNormal{\sp}^{r-2} \sp + \nabla \pi = \mu_\phi \nabla \phi + \mu_\psi \nabla \psi,
\quad \text{in}\ \ L^{r'}_{\mathrm{loc}} \BracketBig{0,\infty;\bm{L}^{r'}(\Omega)}.
\]

Finally, we prove the global uniform boundedness of $(\sp,\phi,\psi,\mu_\phi,\mu_\psi)$ under the additional assumption \eqref{additional_requirements_uniform_boundedness}.
It follows from the energy identity \eqref{energyid} and Young's inequality that
\begin{align*}
    & \frac{\mathrm{d}}{\mathrm{d} t} E_{\mathrm{tot}} \BracketBig{ \sp(t),\phi(t),\psi(t) } + \int_{\Omega} \big(\nu (\phi(t),\psi(t)) |\sp(t)|^2 + \eta (\phi(t),\psi(t)) |\sp(t)|^r\big) \, \mathrm{d} x
    \nonumber \\
    & \qquad + \int_{\Omega} \big( m_{\phi}(\phi(t)) |\nabla \mu_\phi(t)|^2 +m_{\psi}(\psi(t)) |\nabla \mu_\psi(t)|^2\big) \, \mathrm{d} x \nonumber \\
    & \quad \le \AbsBig{ \overline{\phi} - c } \NormNormal{\sigma_1(\phi)}_{L^{\frac{6}{5}}(\Omega)} \NormNormal{\mu_\phi}_{L^6(\Omega)} \nonumber \\
    & \quad \le C \AbsBig{ \overline{\phi} - c } \NormNormal{\sigma_1(\phi)}_{L^{\frac{6}{5}}(\Omega)} \BracketBig{ 1 + \NormNormal{\nabla \mu_\phi} } \nonumber \\
    & \quad \le \frac{1}{2} \underline{m_\phi} \NormNormal{\nabla \mu_\phi}^2 + \AbsBig{ \overline{\phi} - c } \NormNormal{\sigma_1(\phi)}_{L^{\frac{6}{5}}(\Omega)} \BracketMax{ C + \frac{C^2}{2 \underline{m_\phi}} \AbsBig{ \overline{\phi} - c } \NormNormal{\sigma_1(\phi)}_{L^{\frac{6}{5}}(\Omega)} } \\
    & \quad \le \frac{1}{2} \underline{m_\phi} \NormNormal{\nabla \mu_\phi}^2 + C_{\text{diss}} \AbsBig{ \overline{\phi} - c } \NormNormal{\sigma_1(\phi)}_{L^{\frac{6}{5}}(\Omega)},
    \quad \text{for a.a}\ t>0,
\end{align*}
that is,
\begin{align}
    & \frac{\mathrm{d}}{\mathrm{d} t} E_{\mathrm{tot}} \BracketBig{ \sp(t),\phi(t),\psi(t) } + \int_{\Omega} \big(\nu (\phi(t),\psi(t)) |\sp(t)|^2 + \eta (\phi(t),\psi(t)) |\sp(t)|^r\big) \, \mathrm{d} x
    \nonumber \\
    & \qquad + \int_{\Omega} \left(\frac{1}{2} m_{\phi}(\phi(t)) |\nabla \mu_\phi(t)|^2 +m_{\psi}(\psi(t)) |\nabla \mu_\psi(t)|^2\right)\, \mathrm{d} x \nonumber \\
    & \quad \le C_{\text{diss}} \AbsBig{ \overline{\phi} - c } \NormNormal{\sigma_1(\phi)}_{L^{\frac{6}{5}}(\Omega)}, \quad \text{for a.a}\ t>0.
    \label{energy_decreasing_dif}
\end{align}
Here, we have used the boundedness of $\sigma_1$ and the positive constant $C_{\text{diss}}$ is given by
\[
C_{\text{diss}} \stackrel{\mathrm{def}}{=} C + \frac{C^2}{\underline{m_\phi}} \AbsNormal{\Omega}^{\frac{5}{6}} \sigma_1^*.
\]
With the aid of \eqref{energy_decreasing_dif}, we see that the $L^1\BracketBig{0,\infty;L^{\frac{6}{5}}(\Omega)}$-boundedness of $\BracketBig{\overline{\phi}-c} \sigma_1(\phi)$ implies the desired global-in-time estimates.

The proof of Theorem \ref{ws_chdf_thm} is complete.
\qed

\section{Convergence to a Single Equilibrium}
\label{proof_convergence_equilibrium}
The proof of Theorem \ref{ws_chdf_convergence_singleton} is based on some preliminary results. First of all, we need to analyze the $\omega$-limit set of a given global weak solution.
\bl \label{property_equilibrium}
Let the assumptions of Theorem \ref{ws_chdf_thm} hold and $(\sp,\phi,\psi,\mu_\phi,\mu_\psi)$ be a global weak solution to the problem \eqref{eq:chdf1}--\eqref{eq:chdfi} in the sense of Definition \ref{ws_chdf_def}, satisfying the integrability condition \eqref{additional_requirements_uniform_boundedness}.

(1) For any $T>0$, it holds
\begin{align}
\sp(t+\cdot) \to  \bm{0} \quad \text{strongly in } L^r\BracketBig{0,T;\bm{L}^r(\Omega)} \ \text{ as }t\to\infty,\label{translation}
\end{align}

(2) If $\alpha=0$,
$\omega(\phi,\psi) \subset \bm{\mathcal{S}}$ is bounded in $W \times W$;\\ if $\alpha>0$, $\omega(\sp,\phi,\psi) \subset \SetNormal{\bm{0}} \times \bm{\mathcal{S}}$ is bounded in
$\SetNormal{\bm{0}} \times W \times W$, and, in addition,
\begin{align}
    \sp(t)\to \mathbf 0,\quad\text{strongly  in }\bm{L}^2_\sigma(\Omega)\text{ as }t\to\infty.
\end{align}
For $\alpha\geq 0$, there exists $\delta_1 \in (0,\frac14)$ such that
\begin{equation} \label{separation_equilibrium_lemma}
    \NormNormal{\phi_\infty}_{L^\infty(\Omega)} \le 1 - 2\delta_1, \quad  \NormBig{\psi_\infty - \frac{1}{2}}_{L^\infty(\Omega)} \le \frac{1}{2} - 2\delta_1,
\end{equation}
for all $(\phi_\infty,\psi_\infty) \in \omega(\phi,\psi)$ (resp. $(\bm{0},\phi_\infty,\psi_\infty) \in \omega(\sp,\phi,\psi)$). \smallskip

(3) The trajectories $(\phi(t),\psi(t))$ are precompact in $V \times V$. In addition,
\begin{itemize}
    \item if $\alpha=0$, it holds
        \begin{align*}
            \omega(\phi,\psi) = \SetNormal{ (\widetilde{\phi},\widetilde{\psi}) \in (V \cap L^\infty(\Omega))^2 : \exists\, t_n \to \infty \ \text{such that $\phi(t_n) \to \widetilde{\phi}$ and $\psi(t_n) \to \widetilde{\psi}$ in $V$} },
        \end{align*}
        and
        \begin{align}
            \lim\limits_{t \to \infty} \mathrm{dist}_{V \times V} \BracketBig{ (\phi(t),\psi(t)),\omega(\phi,\psi) } =0; \label{dist_alpha=0-b}
        \end{align}
    \item if $\alpha >0$, it holds
        \begin{align*}
            \omega(\sp,\phi,\psi) = \SetNormal{ (\bm{0},\widetilde{\phi},\widetilde{\psi}) \in \{ \bm{0} \} \times (V \cap L^\infty(\Omega))^2 & : \exists\, t_n \to \infty \ \text{such that $\phi(t_n) \to \widetilde{\phi}$ and} \\
            & \quad \text{$\psi(t_n) \to \widetilde{\psi}$ in $V$, $\sp(t_n) \to \bm{0}$ in $\bm{L}^2_\sigma(\Omega)$} },
        \end{align*}
        and
        \begin{align}
        \lim\limits_{t \to \infty} \mathrm{dist}_{\bm{L}^2_\sigma(\Omega) \times V \times V} \BracketBig{ (\sp(t),\phi(t),\psi(t)),\omega(\sp,\phi,\psi) } =0. \label{dist_alpha>0-b}
        \end{align}
\end{itemize}
\el
\noindent \textbf{Proof.}
Given $\BracketBig{\widetilde{\phi},\widetilde{\psi}} \in \omega(\phi,\psi)$ if $\alpha=0$ (resp. $\BracketBig{\widetilde {\sp},\widetilde{\phi},\widetilde{\psi}} \in \omega(\sp,\phi,\psi)$ if $\alpha>0$), we consider an increasing sequence $t_n \to \infty$ such that $(\phi(t_n),\psi(t_n))$ converges to $(\widetilde{\phi},\widetilde{\psi})$
weakly in $V \times V$ and
strongly in $H \times H$, and if $\alpha>0$, $\sp(t_n)$ converges to $\widetilde{\sp}$ weakly in $\bm{L}^2_\sigma(\Omega)$.
For brevity, we define the sequence of time-shifted solutions as
\[
(\sp_n,\phi_n,\psi_n,\mu_{\phi,n},\mu_{\psi,n}) \stackrel{\mathrm{def}}{=} (\sp(t+t_n),\phi(t+t_n),\psi(t+t_n),\mu_\phi(t+t_n),\mu_\psi(t+t_n)),
\]
which solve the following system
\begin{numcases}{}
    - \alpha \InnerNormal{\sp_n}{\partial_t \bm{v}} + \PairingNormal{ \eta (\phi_n,\psi_n) \AbsNormal{\sp_n}^{r-2} \sp_n }{ \bm{v} }_{ L^{r'} \BracketNormal{ 0,T;\bm{L}^{r'}(\Omega) },L^r \BracketNormal{ 0,T;\bm{L}^r(\Omega) } }  \nonumber \\
    \quad + \InnerNormal{\nu (\phi_n,\psi_n) \sp_n}{\bm{v}} = - \InnerNormal{ \phi_n \nabla \mu_{\phi,n} + \psi_n \nabla \mu_{\psi,n} }{ \bm{v} } \label{eq:df1n}, & $\forall\,\bm{v} \in C^\infty_0 \BracketBig{(0,T);\bm{L}^r_\sigma(\Omega)}$,\\
    \PairingNormal{ \partial_t \phi_n + \sigma_1(\phi_n) \BracketBig{\overline{\phi_n} -c} }{ \zeta }_{L^2 (0,T;V^*),L^2(0,T;V)} - \InnerNormal{ \phi_n \sp_n }{\nabla \zeta}_{Q_T} \nonumber \\
    \quad = - \InnerNormal{ m_{\phi} (\phi_n) \nabla \mu_{\phi,n} }{ \nabla \zeta }_{Q_T} \label{eq:ch1n}, & $\forall\, \zeta \in L^2\BracketBig{0,T;V}$, \\
    \mu_{\phi,n} = -\Delta \phi_n + \sigma_2 \mathcal{N} \BracketBig{\phi_n - \overline{\phi_n}} + F_{\phi}'(\phi_n) + \partial_{\phi} G(\phi_n,\psi_n) \label{eq:ch2n}, & a.e. in $Q_T$,\\
    \PairingNormal{ \partial_t \psi_n }{ \zeta }_{L^2 (0,T;V^*),L^2(0,T;V)} - \InnerNormal{ \psi_n \sp_n }{\nabla \zeta}_{Q_T} \nonumber \\
    \quad = - \InnerNormal{ m_{\psi} (\psi_n) \nabla \mu_{\psi,n} }{ \nabla \zeta }_{Q_T} \label{eq:ch3n}, & $\forall\, \zeta \in L^2\BracketBig{0,T;V}$, \\
    \mu_{\psi,n} = - \beta \Delta \psi_n + F_{\psi}'(\psi_n) + \partial_{\psi} G(\phi_n,\psi_n) \label{eq:ch4n}, & a.e. in $Q_T$,
\end{numcases}
for any $T>0$.
Define
\[
e(t) \stackrel{\mathrm{def}}{=} E_{\mathrm{tot}}(\sp(t),\phi(t),\psi(t))
+ C_{\text{diss}} R_{\text{diss}}(t), \quad R_{\text{diss}}(t) \stackrel{\mathrm{def}}{=} \int_t^\infty \AbsBig{\overline{\phi} - c} \NormNormal{\sigma_1(\phi)}_{L^{\frac{6}{5}}(\Omega)} \ \mathrm{d} \tau,
\]
where $C_{\text{diss}}$ is the positive constant appearing in
\eqref{energy_decreasing_dif}.
We infer from \eqref{energy_decreasing_dif} that
\begin{align}
    & e(t) + \int_{Q_{(s,t)}} \left( \nu (\phi,\psi) |\sp|^2 + \eta (\phi,\psi) |\sp|^r + \frac{1}{2} m_{\phi}(\phi) |\nabla \mu_\phi|^2 +m_{\psi}(\psi) |\nabla \mu_\psi|^2\right) \, \mathrm{d} (x,\tau) \le e(s), \label{energy_decreasing}
\end{align}
for all $t \in [0,\infty)$ and almost all $s\in[0,t)$ including $s=0$. Recalling the definition of $R_{\text{diss}}$, we see that $\lim\limits_{t \to \infty} R_{\text{diss}}(t) = 0$. This, combined with the fact that the function $e$ is nonincreasing in time and is bounded from below, guarantees that
\begin{align}
\lim_{t\to\infty} E_{\mathrm{tot}}(\sp(t),\phi(t),\psi(t))=E_\infty,
\label{convr}
\end{align}
for some $E_\infty\in \mathbb{R}$.

From the energy inequality \eqref{energy_decreasing}, we see that there exists
$(\sp^*,\phi^*,\psi^*,\mu_\phi^*,\mu_\psi^*)
$ such that, for any given $T>0$, it holds
\begin{align*}
    & \phi_n \rightharpoonup \phi^*,\ \  \psi_n \rightharpoonup \psi^* \quad \text{weakly in $L^2(0,T;W) \cap H^1(0,T;V^*)$}, \\
    & \phi_n \stackrel{*}{\rightharpoonup} \phi^*,\ \  \psi_n \stackrel{*}{\rightharpoonup} \psi^* \quad \text{weakly-* in $L^\infty(0,T;V)$}, \\
    & \phi_n \to \phi^*,\ \  \psi_n \to \psi^* \quad \text{in $L^2(0,T;H^{s+1}(\Omega)) \cap C([0,T];H^s(\Omega))$,\ \  $\forall\, s \in [0,1)$}, \\
    & \mu_{\phi,n} \rightharpoonup \mu_\phi^*,\ \  \mu_{\psi,n} \rightharpoonup \mu_\psi^* \quad \text{weakly in $L^2(0,T;V)$}, \\
    & \sp_n \rightharpoonup \sp^* \ \ \text{weakly in $L^r(0,T;\bm{L}^r_\sigma(\Omega))$, \ \ \ and }\  \alpha \sp_n \rightharpoonup \alpha \sp^* \ \  \text{weakly in $W^{1,r'}(0,T;\bm{L}^{r'}_\sigma(\Omega))$}.
\end{align*}
With a standard argument, we can show that $(\sp^*,\phi^*,\psi^*,\mu_\phi^*,\mu_\psi^*)$ satisfies, for any fixed $T>0$,
\begin{numcases}{}
    \PairingNormal{ \partial_t \phi^* + \sigma_1(\phi^*) \BracketBig{\overline{\phi^*} -c} }{ \zeta }_{L^2 (0,T;V^*),L^2(0,T;V)} - \InnerNormal{ \phi^* \sp^* }{\nabla \zeta}_{Q_T} \nonumber \\
    \quad = - \InnerNormal{ m_{\phi} (\phi^*) \nabla \mu_\phi^* }{ \nabla \zeta }_{Q_T} \label{eq:ch1limit}, & $\forall\, \zeta \in L^2\BracketBig{0,T;V}$, \\
    \mu_{\phi}^* = -\Delta \phi^* + \sigma_2 \mathcal{N} \BracketBig{\phi^* - \overline{\phi^*}} + F_{\phi}'(\phi^*) + \partial_{\phi} G(\phi^*,\psi^*) \label{eq:ch2limit}, & a.e. in $Q_T$,\\
    \PairingNormal{ \partial_t \psi^* }{ \zeta }_{L^2 (0,T;V^*),L^2(0,T;V)} - \InnerNormal{ \psi^* \sp^* }{\nabla \zeta}_{Q_T} \\
    \quad = - \InnerNormal{ m_{\psi} (\psi^*) \nabla \mu_\psi^* }{ \nabla \zeta }_{Q_T} \label{eq:ch3limit}, & $\forall\, \zeta \in L^2\BracketBig{0,T;V}$, \\
    \mu_{\psi}^* = - \beta \Delta \psi^* + F_{\psi}'(\psi^*) + \partial_{\psi} G(\phi^*,\psi^*) \label{eq:ch4limit}, & a.e. in $Q_T$,\\
    \phi^*(0)=\widetilde\phi,\quad  \psi^*(0)={\widetilde \psi}, & a.e. in $\Omega$.
\end{numcases}
From \eqref{energy_decreasing} and the convergence of $E_{\mathrm{tot}}(\sp,\phi,\psi)$, we infer that
\begin{align}
    & E_{\infty} + \int_{Q_{(s,t)}} \nu (\phi^*,\psi^*) |\sp^*|^2 + \eta (\phi^*,\psi^*) |\sp^*|^r + \frac{1}{2} m_{\phi}(\phi^*) |\nabla \mu_\phi^*|^2 +m_{\psi}(\psi^*) |\nabla \mu_\psi^*|^2 \ \mathrm{d} (x,\tau) \le E_\infty , 
    \notag
\end{align}
for all $t \in [0,\infty)$ and almost all $s\in[0,t)$ including $s=0$. Since $s,t$ are arbitrary, we can conclude
\begin{align*}
    & \nabla \mu_{\phi}^*(t)\equiv \mathbf{0},\quad \nabla \mu_{\psi}^*(t)\equiv\mathbf 0,\quad \sp^*(t)\equiv \mathbf{0}, \quad \text{for a.a.}\ t>0.
\end{align*}
Moreover, we can conclude the convergence  \eqref{translation} for the time-translation $\sp(t+\cdot)$. Indeed, from \eqref{energy_decreasing} we also obtain that $\int_0^T\NormBig{\sp(t+\tau)}^r_{\mathbf \bm{L}^r(\Omega)}d \tau\to0$ as $t\to\infty$, thanks to the convergence \eqref{convr}.
Moreover, it follows from the almost everywhere pointwise convergence of $\phi$ and
\[ \int_0^T \NormNormal{\sigma_1(\phi_n)}_{L^{\frac{6}{5}}(\Omega)} \AbsBig{\overline{\phi_n}-c} \ \mathrm{d} \tau \to 0, \]
which is due to the assumption \eqref{additional_requirements_uniform_boundedness},
that
\[ \NormNormal{ \sigma_1(\phi^*) \BracketBig{\overline{\phi^*}-c} }_{L^2(0,T;L^{\frac{6}{5}}(\Omega))}=0,
\quad \forall\, T>0.
\]
If $\alpha>0$, since $\sp^*\in W^{1,r'}(0,T;\bm{L}^{r'}(\Omega))\hookrightarrow C([0,T];\bm{L}^{r'}(\Omega))$, we can conclude $\widetilde \sp=\sp^*(0)=\mathbf 0$ almost everywhere in $\Omega$.

Therefore, the system \eqref{eq:ch1limit}--\eqref{eq:ch4limit} reduces to
\begin{align*}
    & \partial_t \phi^* = \partial_t \psi^* = 0,\quad  \text{ in }V' \ \text{ for a.a. }t>0, \\
    & \mu_{\phi,\infty} = -\Delta \phi^* + \sigma_2 \mathcal{N} \BracketBig{\phi^* - \overline{\phi^*}} + F_{\phi}'(\phi^*) + \partial_{\phi} G(\phi^*,\psi^*), \quad \text{a.e. in $Q_T$}, \\
    & \mu_{\psi,\infty} = - \beta \Delta \psi^* + F_{\psi}'(\psi^*) + \partial_{\psi} G(\phi^*,\psi^*), \quad \text{a.e. in $Q_T$}, \\
    & \partial_{\bm{n}} \phi^* = \partial_{\bm{n}} \psi^* = 0, \quad \text{a.e. on $S_T$},\\&
    \phi^*(0)=\widetilde\phi,\quad
    \psi^*(0)=\widetilde \psi,\quad\text{ a.e. in }\Omega,
\end{align*}
for some constants $\mu_{\phi,\infty},\mu_{\psi,\infty} \in \mathbb{R}$ with $(\mu_\phi^*,\mu_\psi^*) \equiv (\mu_{\phi,\infty},\mu_{\psi,\infty})$.
This implies that $(\widetilde \phi,\widetilde\psi) \in \bm{\mathcal{S}}$ with $(\phi^*(t),\psi^*(t)) \equiv (\widetilde \phi,\widetilde\psi)$ almost everywhere in $\Omega$ for all $t \ge 0$. Therefore, we have   $\omega(\phi,\psi)\subset \mathcal S$ for $\alpha=0$, and $\omega(\sp,\phi,\psi)\subset \{\mathbf 0 \}\times \mathcal S$ for $\alpha>0$.
Moreover, as $n\to \infty$,
\begin{align}
    & \mathop{\mathrm{sup}}\limits_{t \in [0,T]} \NormNormal{ \phi_n(t) - \widetilde\phi } + \NormNormal{ \psi_n(t) - \widetilde\psi } \to 0,  \\
    & \label{phic}(\phi_n(t),\psi_n(t)) \to (\widetilde\phi,\widetilde\psi) \ \text{strongly in $V \times V$ for a.a. $t>0$}.
\end{align}
When $\alpha=0$, exploiting the convergence $E_{\text{tot}}(\sp(t),\phi(t),\psi(t))\to E_\infty$ as $t\to\infty$ (there is no dependence of $\sp(t)$ since $\alpha=0$), and the above convergence results, we find, for any $s\geq 0$,
\[
\frac{1}{2}  \NormNormal{\nabla \phi(s+t_n)}^2 + \frac{\beta}{2}\NormNormal{\nabla \psi(s+t_n)}^2 \to E_\infty - \int_\Omega \left( F_\phi(\widetilde\phi) + F_\psi(\widetilde\psi) + G(\widetilde\phi,\widetilde\psi) \right)\, \mathrm{d}x,\quad \text{as}\ n\to \infty,
\]\text{which entails, together with \eqref{phic}}, for almost any $s>0$,
\[ \frac{1}{2}  \NormNormal{\nabla \widetilde\phi}^2 + \frac{\beta}{2}\NormNormal{\nabla \widetilde\psi}^2  =
\lim_{n\to \infty} \left(\frac{1}{2}  \NormNormal{\nabla \phi_n(s)}^2 + \frac{\beta}{2}\NormNormal{\nabla \psi_n(s)}^2\right)
= E_\infty - \int_\Omega \left( F_\phi(\widetilde\phi) + F_\psi(\widetilde\psi) + G(\widetilde\phi,\widetilde\psi) \right)\, \mathrm{d}x,
\] which implies $\frac{1}{2}  \NormNormal{\nabla \phi(t_n)}^2 + \frac{\beta}{2}\NormNormal{\nabla \psi(t_n)}^2 \to \frac{1}{2}  \NormNormal{\nabla \widetilde\phi}^2 + \frac{\beta}{2}\NormNormal{\nabla \widetilde\psi}^2$, so that
\[ (\phi(t_n),\psi(t_n)) \to (\widetilde\phi,\widetilde \psi), \quad \text{strongly in $V \times V$}. \]
Here, we have also used the convergence (resp. conservation) of the mass $\overline{\phi}(t_n)$ (resp. $\overline{\psi}(t_n)$) and Poincar\'e--Wirtinger inequality. The above result yields the precompactness of the trajectories $(\phi, \psi)$ in $V\times V$. When $\alpha>0$, we should also prove $\sp \to \bm{0}$ in $\bm{L}^2_\sigma(\Omega)$. It follows from \eqref{testdf} that
\begin{align*}
    & \alpha \NormNormal{\sp(t)}^2 + \int_{Q_{(s,t)}} \big(\nu(\phi,\psi) \AbsNormal{\sp}^2 + \eta(\phi,\psi) \AbsNormal{\sp}^r\big)\, \mathrm{d} (x,\tau) \\
    & \quad = \alpha \NormNormal{\sp(s)}^2 - \InnerNormal{\phi \nabla \mu_\phi + \psi \nabla \mu_\psi}{\sp}_{Q_{s,t}} \\
    & \quad \le \alpha \NormNormal{\sp(s)}^2 + \BracketBig{ \NormNormal{\nabla \mu_\phi}_{L^2(Q_{(s,\infty)})} + \NormNormal{\nabla \mu_\psi}_{L^2(Q_{(s,\infty)})} } \NormNormal{\sp}_{L^2(Q_{(s,\infty)})}.
\end{align*}
By the Chebyshev inequality and the $L^2(Q)$-integrability of $\sp$, we see that for any $\epsilon \in (0,1)$, there exists $T_\epsilon>0$ such that
\[
\NormNormal{\sp(T_\epsilon)} \le \epsilon \quad  \mathrm{and} \quad  \NormNormal{\sp}_{L^2(Q_{(T_\epsilon,\infty)})} \le \epsilon^2.
\]
This implies
\[
\NormNormal{\sp(t)} \le C \epsilon, \quad \forall\,t\ge T_\epsilon,
\]
namely, $\|\sp(t)\|\to 0$ as $t\to\infty$. In addition, for $\alpha\geq 0$, we can apply the contradiction argument in \cite[Section 4]{GP2025} to show the convergence \eqref{dist_alpha=0-b} (resp. \eqref{dist_alpha>0-b}).

The uniform boundedness of $\omega(\sp,\phi,\psi)$ in $\{\mathbf 0\}\times W \times W$ (resp. $\omega(\phi,\psi)$ in $ W \times W$) is based on three facts: the energy inequality yields uniform $H^1(\Omega)$-estimates for every $(\widetilde\phi, \widetilde\psi)$ in the $\omega$-limit set, the corresponding chemical potentials $\mu_{\widetilde\phi}$, $\mu_{\widetilde\psi}$ are constants, and the convexity of $F_\psi$, $F_\psi$. The details are similar to those of \cite[Section 4]{GP2025} and will not be repeated here.
The strict separation property of every element in the $\omega$-limit set is a direct consequence of \cite[Lemma A.1]{CG2021} and the singularity of $F_\phi'$, $F_\psi'$. Furthermore, with the above mentioned uniform $W$-bounds, we can show the uniform separation property \eqref{separation_equilibrium_lemma}, following the argument for \cite[Lemma 3.11]{GGPS23}.
%
The proof is complete.
\qed
\medskip

Then, we establish the validity of the strict separation property for $\phi$ and $\psi$ in the set of ``good times''.

\bl \label{longtime_separation}
Let the assumptions of Theorem \ref{ws_chdf_thm} hold and $(\sp,\phi,\psi,\mu_\phi,\mu_\psi)$ be a global weak solution to problem \eqref{eq:chdf1}--\eqref{eq:chdfi} in the sense of Definition \ref{ws_chdf_def}, satisfying the integrability condition \eqref{additional_requirements_uniform_boundedness}.

(1) Consider the set
\[ A_\delta(t) \stackrel{\mathrm{def}}{=} \SetBig{x\in \Omega: \ \AbsNormal{\phi(x,t)}\ge 1- \frac{3}{2} \delta_1 \ \text{ or } \ \AbsNormal{\psi(x,t)-\frac{1}{2}}\ge \frac{1}{2}- \frac{3}{2} \delta_1}, \quad t\ge 0, \]
where $\delta_1 > 0$ is given as in \eqref{separation_equilibrium_lemma}. Then we have
\begin{align}
    \lim\limits_{t \to \infty} \AbsNormal{A_\delta(t)} = 0. \label{weakl}
\end{align}

(2) For any $M>0$, there exist $\delta \in (0, \delta_1)$ and $T_M > 0$ such that
\begin{equation} \label{longtime_separation_good_time}
    \sup\limits_{t \in A_{M}(T_M)} \NormNormal{\phi(t)}_{L^\infty(\Omega)} \le 1 - \delta, \quad \sup\limits_{t \in A_{M}(T_M)} \NormBig{\psi(t) - \frac{1}{2}}_{L^\infty(\Omega)} \le \frac{1}{2} - \delta,
\end{equation}
where
\[ A_M(T) \stackrel{\mathrm{def}}{=} \SetNormal{ t\ge T : \ \NormNormal{\nabla \mu_\phi(t)}^2 +\NormNormal{\nabla \mu_\psi(t)}^2 \le M^2 }\]
is the (measurable) set of ``good times''.
\el

\noindent \textbf{Proof.}
We adopt the strategy of \cite[Lemma 3.4]{GP2025} and focus mainly on the differences in the argument. Modifications are necessary due to the coupling structure of the Cahn--Hilliard system and the nonlocal interaction term.

Concerning the asymptotic separation property \eqref{weakl}, the main differences here are
due to the presence of two phase-field variables. To this end, for any $(\phi_\infty,\psi_\infty) \in \omega(\phi,\psi)$ (resp. $(\bm{0},\phi_\infty,\psi_\infty) \in \omega(\sp,\phi,\psi)$), we introduce the following set:
\begin{align*}
    & B_\delta^{\phi_\infty,\psi_\infty}(t) \stackrel{\mathrm{def}}{=} \SetNormal{x\in \Omega: \ \AbsNormal{\phi(x,t) - \phi_\infty(x)} \ge \frac{1}{2} \delta_1 \ \text{or} \ \AbsNormal{\psi(x,t)-\psi_\infty(x)} \ge \frac{1}{2} \delta_1}.
\end{align*}
It follows from \eqref{separation_equilibrium_lemma} that, for all $t \ge 0$ and almost all $x \in A_\delta(t)$, we have either
\[ 1 - \frac{3}{2} \delta_1 \le \AbsNormal{\phi(x,t)} \le \AbsNormal{\phi_\infty(x)} + \AbsNormal{\phi(x,t)-\phi_\infty(x)} \le 1 - 2 \delta_1 + \AbsNormal{\phi(x,t) - \phi_\infty(x)}, \]
or
\[ \frac{1}{2} - \frac{3}{2} \delta_1 \le \AbsNormal{\psi(x,t) - \frac{1}{2}} \le \AbsNormal{\psi_\infty(x)-\frac{1}{2}} + \AbsNormal{\psi(x,t)-\psi_\infty(x)} \le \frac{1}{2} - 2 \delta_1 + \AbsNormal{\psi(x,t) - \psi_\infty(x)}. \]
These inequalities imply that $\AbsNormal{\phi(x,t) - \phi_\infty(x)} \ge \frac{1}{2} \delta_1$ or $\AbsNormal{\psi(x,t) - \psi_\infty(x)} \ge \frac{1}{2} \delta_1$, which entail
\[ A_\delta (t) \subset B_\delta^{\phi_\infty,\psi_\infty}(t), \quad \forall\, t \ge 0. \]
Using the Chebyshev inequality, we obtain
\begin{equation*}
    \AbsNormal{A_\delta(t)} \le \AbsNormal{B_\delta^{\phi_\infty,\psi_\infty}(t)} \le \frac{4}{\delta_1^2} \BracketBig{ \NormNormal{\phi(t) - \phi_\infty}^2 + \NormNormal{\psi(t) - \psi_\infty}^2 },\quad \forall\, t\geq 0.
\end{equation*}
Taking the infimum over all elements $(\phi_\infty,\psi_\infty)$ in the $\omega$-limit set yields
\begin{align*}
    & \AbsNormal{A_\delta(t)} \le \frac{4}{\delta_1^2} \mathrm{dist}^2_{H \times H} \BracketBig{(\phi(t),\psi(t)),\omega(\phi,\psi)} \to 0, \quad \text{if $\alpha=0$}; \\
    & \AbsNormal{A_\delta(t)} \le \frac{4}{\delta_1^2} \mathrm{dist}^2_{\bm{L}^2_\sigma(\Omega) \times H \times H} \BracketBig{(\sp(t),\phi(t),\psi(t)),\omega(\sp,\phi,\psi)} \to 0, \quad \text{if $\alpha>0$},
\end{align*}
as $t \to \infty$. This allows us to conclude \eqref{weakl}.

For the second result concerning the uniformly strict separation property on the set of ``good times'', we apply a De Giorgi's iteration scheme (see \cite{GalP}). As in \cite{GP2025}, we define
\[k_n = 1 - \delta - \frac{\delta}{2^n}, \quad \forall\, n\ge0,\]
for some fixed $\delta \in (0, \delta_1)$, and for $t\geq 0$
\[ \phi_n(x,t) \stackrel{\mathrm{def}}{=} (\phi-k_n)^+, \quad A_n(t) \stackrel{\mathrm{def}}{=} \SetNormal{x\in \Omega : \ \phi(x,t)\ge k_n}, \quad y_n(t) = \int_{A_n(t)} 1 \ \mathrm{d} x. \]
Then, by definition, $0\leq \phi_n\leq 2\delta$.
Testing equation \eqref{eq:chdf4} with $\phi_n-\overline{\phi_n}$, we get
\[ \NormNormal{\nabla \phi_n}^2 + \int_\Omega F_\phi' (\phi) (\phi_n - \overline{\phi_n}) \ \mathrm{d} x = - \InnerNormal{\partial_\phi G (\phi,\psi) + \sigma_2 \mathcal{N} (\phi-\overline{\phi})}{\phi_n - \overline{\phi_n}} + \InnerNormal{\mu_\phi}{\phi_n - \overline{\phi_n}}.
\]
Note that the second term on the right-hand side involving the coupling and nonlocal interaction can be controlled by
\begin{align*}
     \AbsMax{ \InnerBig{\partial_\phi G (\phi,\psi) + \sigma_2 \mathcal{N} (\phi-\overline{\phi})}{\phi_n - \overline{\phi_n}} }& = \AbsMax{ \InnerBig{ \partial_\phi G (\phi,\psi) - \overline{\partial_\phi G (\phi,\psi)} + \sigma_2 \mathcal{N} (\phi-\overline{\phi})}{\phi_n} } \\
    & \quad \le C \int_\Omega \phi_n \ \mathrm{d} x
    \le 2C \delta |\Omega|^\frac16 y_n(t)^{\frac{5}{6}}.
\end{align*}
The other terms can be estimated by the same argument as in \cite[Section 5]{GP2025}. Hence, we can obtain
\[
y_{n+1}(t) \le \frac{1}{\delta^3} 2^{3n+3} C(\delta,M)^{\frac{9}{10}} y_n(t)^{\frac{29}{20}},
\qquad \forall\, n \ge 0,
\]
for any $t\in A_M(T)$, with $T\geq \widetilde{T}$ for $\widetilde{T}$ sufficiently large.
Using \eqref{weakl} together with the well-known geometric lemma (see \cite[Lemma~A.2]{GP2025}), we eventually obtain $\sup_{t\in A_M(T)}\|(\phi(t)-(1-\delta))^+\|_{L^\infty(\Omega)}=0$. Analogously, we can get $\sup_{t\in A_M(T)}\|(\phi(t)-(-1+\delta))^+\|_{L^\infty(\Omega)}=0$. In addition, the same procedure applies to $\psi$. This completes the proof of \eqref{longtime_separation_good_time}.
\qed
\smallskip

The last crucial tool is the following extended version of the {\L}ojasiewicz--Simon inequality for coupled elliptic system \eqref{ch_stationary_phi}--\eqref{ch_stationary_boundary} (see \cite{O2025} for a proof).
\bl[{\L}ojasiewicz--Simon inequality] \label{ls_inequality}
Let $\Omega \subset \mathbb{R}^d$, $d\in\{2,3\}$, be a bounded domain with a smooth boundary $\partial\Omega$.
Suppose that the assumptions (\textbf{H1}) and (\textbf{H1*}) are satisfied.
Given a global weak solution $(\sp,\phi,\psi,\mu_\phi,\mu_\psi)$ to problem \eqref{eq:chdf1}--\eqref{eq:chdfi} in the sense of Definition \ref{ws_chdf_def}, satisfying the integrability condition \eqref{additional_requirements_uniform_boundedness}, for any $(\phi_\infty,\psi_\infty) \in \omega(\phi,\psi)$ if $\alpha=0$ (resp. $(\bm{0},\phi_\infty,\psi_\infty) \in \omega(\sp,\phi,\psi)$ if $\alpha>0$),
there exist constants $\theta \in (0,\frac12)$ and $C,\varpi > 0$ such that for $\widetilde{\phi},\widetilde{\psi} \in V\cap L^\infty(\Omega)$ satisfying
\[
\NormNormal{\widetilde{\phi}-\phi_\infty}_{V} \leq \varpi, \quad \NormNormal{\widetilde{\psi}-\psi_\infty}_{V} \leq \varpi,\quad
\NormNormal{\widetilde{\phi}}_{L^{\infty}(\Omega)} \leq 1 -\frac{3}{2}\delta_1,\quad  \NormBig{\widetilde{\psi}-\frac{1}{2}}_{L^{\infty}(\Omega)} \leq \frac{1}{2} -\frac{3}{2} \delta_1,\]
where $\delta_1\in (0,1)$ is the constant determined in Lemma \ref{property_equilibrium},
the following inequality holds
\begin{align}
   &  \big|E_{\mathrm{free}}(\widetilde{\phi},\widetilde{\psi})-E_{\mathrm{free}}(\phi_\infty,\psi_\infty)\big|^{1-\theta}\notag \\
   &\quad  \le C \left( \|\widetilde{\mu}_\phi-\overline{\widetilde{\mu}_\phi}\|_{V^*_0}+ \|\widetilde{\mu}_\psi-\overline{\widetilde{\mu}_\psi}\|_{V^*_0} +  \left(|\overline{\widetilde{\phi}}-\overline{\phi_\infty}| + |\overline{\widetilde{\psi}}-\overline{\psi_\infty}|\right)^{1-\theta}\right), \label{eq:ls_inequality}
\end{align}
with $\widetilde{\mu}_\phi=-\Delta \widetilde{\phi} + F_\phi'(\widetilde{\phi}) + \sigma_2 \mathcal{N}(\widetilde{\phi} -\overline{\widetilde{\phi}}) + \partial_\phi G(\widetilde{\phi},\widetilde{\psi})$
and $\widetilde{\mu}_\psi=-\beta \Delta \widetilde{\psi} + F_\psi'(\widetilde{\psi}) + \partial_\psi G(\widetilde{\phi},\widetilde{\psi})$.
\el

We are ready to prove Theorem \ref{ws_chdf_convergence_singleton}, namely, to show the convergence to a single equilibrium for global weak solutions satisfying the decay condition \eqref{spatial_average_decay}, under the assumption that $F_\phi,F_\psi,G$ are real analytic.

\smallskip

\noindent \textbf{Proof of Theorem \ref{ws_chdf_convergence_singleton}.}
Suppose that $(\sp,\phi,\psi,\mu_\phi,\mu_\psi)$ is a global weak solution of
problem \eqref{eq:chdf1}--\eqref{eq:chdfi} satisfying \eqref{spatial_average_decay}.
By the compactness of the $\omega$-limit set, there exists a finite $\delta$-net
$\{(\phi_i,\psi_i)\}_{i=1}^n$ of this set with $\delta = \min\limits_{1 \le i \le n} \SetNormal{\tfrac{\varpi_i}{2}}$,
where $\varpi_i$ is the constant associated with
$(\phi_i,\psi_i)\in\omega(\phi,\psi)$ (resp.
$(\bm{0},\phi_i,\psi_i) \in \omega(\sp,\phi,\psi)$) as given in
Lemma~\ref{ls_inequality},
such that
\begin{align*}
    & \omega(\phi,\psi)\subset \bigcup_{i=1}^n
    B_{V\times V}\!\left((\phi_i,\psi_i),\delta \right),
    \quad \text{if $\alpha=0$}; \\
    & \omega(\sp,\phi,\psi)\subset \SetNormal{\bm{0}} \times
    \bigcup_{i=1}^n
    B_{V\times V}\!\left((\phi_i,\psi_i),\delta \right),
    \quad \text{if $\alpha>0$}.
\end{align*}
Moreover, from Lemma~\ref{property_equilibrium} it follows that there exists
$t_*>0$ such that for $t\geq t_*$,
\begin{align*}
    & \mathrm{dist}_{V \times V} \BracketBig{(\phi(t),\psi(t)),\omega(\phi,\psi)} \le \frac{\delta}{2}, \quad \text{if $\alpha=0$}; \\
    & \mathrm{dist}_{\bm{L}^2_\sigma(\Omega) \times V \times V} \BracketBig{(\sp(t),\phi(t),\psi(t)),\omega(\sp,\phi,\psi)} \le \frac{\delta}{2}, \quad \text{if $\alpha>0$}.
\end{align*}
Let $\theta_i$ be the {\L}ojasiewicz exponent corresponding to $(\phi_i,\psi_i)$ in Lemma~\ref{ls_inequality}, $i \in \mathbb{N}_+ \cap [1,n]$.
We set
\begin{align}
 \zeta = \min \SetMax{\min\limits_{1 \le i \le n}\SetBig{\theta_i},\frac{\rho}{2(1+\rho)}}, \label{constraint}
 \end{align}
where $\rho >0$ is the constant related to the decay rate given in \eqref{spatial_average_decay}.
It follows from \eqref{energy_decreasing}
that
\begin{align*}
    & C_* \BracketMax{
    \int_s^t
    \NormNormal{\nabla \mu_\phi}^2
    + \NormNormal{\nabla \mu_\psi}^2
    + \NormNormal{\sp}_{\bm{L}^2_\sigma(\Omega)}^2
    + \NormNormal{\sp}_{\bm{L}^r_\sigma(\Omega)}^r
    \,\mathrm{d} \tau }
    \le e(s) - e(t) ,
\end{align*}
for all $t \in [0,\infty)$ and almost any $s\in [0,t)$ including $s=0$. By adding to both sides the positive quantity
$C_*\int_s^t (1+\tau)^{-2(1-\zeta)(1+\rho)} \ \mathrm{d} \tau$, we get
\begin{align*}
    & C_* \BracketMax{
    \int_s^t
    \NormNormal{\nabla \mu_\phi}^2
    + \NormNormal{\nabla \mu_\psi}^2
    + \NormNormal{\sp}_{\bm{L}^2_\sigma(\Omega)}^2
    + \NormNormal{\sp}_{\bm{L}^r_\sigma(\Omega)}^r+(1+\tau)^{-2(1-\zeta)(1+\rho)}
    \,\mathrm{d} \tau } \\
    & \quad \le e(s) - e(t) + C_* \BracketNormal{ (1+s)^{-2(1+\rho)(1-\zeta)+1} - (1+t)^{-2(1+\rho)(1-\zeta)+1} }.
\end{align*}
Raising both sides to the power of $2(1-\zeta)$ and then passing to the limit as $t \to \infty$, recalling that $1-\zeta>\frac12$, we obtain
\begin{align*}
    & \BracketMax{
    \int_s^\infty
    \NormNormal{\nabla \mu_\phi}^2
    + \NormNormal{\nabla \mu_\psi}^2
    + \NormNormal{\sp}_{\bm{L}^2_\sigma(\Omega)}^2
    + \NormNormal{\sp}_{\bm{L}^r_\sigma(\Omega)}^r+(1+\tau)^{-2(1-\zeta)(1+\rho)}
    \,\mathrm{d} \tau }^{2(1-\zeta)} \\
    & \quad \le \frac{1}{C_*^{2(1-\zeta)}}
    \left( e(s) - E_\infty + C_* (1+s)^{-(1+\rho)} \right)^{2(1-\zeta)} \\
    & \quad \le
    \frac{4^{2(1-\zeta)}}{C_*^{2(1-\zeta)}}
    \AbsNormal{ E_{\mathrm{free}}(\phi(s),\psi(s)) - E_\infty }^{2(1-\zeta)}
    \bigl(
    \chi_{A_M(T_M)}(s)
    + \chi_{(t_*,\infty)\backslash A_M(T_M)}(s)
    \bigr) \\
    & \qquad
    + \BracketMax{\frac{4}{C_*}}^{2(1-\zeta)}
    \alpha^{2(1-\zeta)} \NormNormal{\sp(s)}^{4(1-\zeta)}
    +
    \BracketMax{\frac{{4}C_{\text{diss}}}{C_*}}^{2(1-\zeta)}
    R_{\text{diss}}^{2(1-\zeta)}(s)\\
    &
    \qquad + {4}^{2(1-\zeta)} (1+s)^{-2(1-\zeta)(1+\rho)}\\
    & \quad \le
    \frac{4^{2(1-\zeta)}}{C_*^{2(1-\zeta)}}
    \AbsNormal{ E_{\mathrm{free}}(\phi(s),\psi(s)) - E_\infty }^{2(1-\zeta)}
    \bigl(
    \chi_{A_M(T_M)}(s)
    + \chi_{(t_*,\infty)\backslash A_M(T_M)}(s)
    \bigr) \\
    & \qquad
    + \BracketMax{\frac{4}{C_*}}^{2(1-\zeta)}
    \alpha (e(0)+C)^{1-2\zeta}
    \NormNormal{\sp(s)}^2 + C (1+s)^{-2(1-\zeta)(1+\rho)},
\end{align*}
for all $s \in (t_*,\infty)$, since $\zeta \in (0,\frac12)$ and $R_{\text{diss}}(s)\le C (1+s)^{-2(1+\rho)}$ (see the second condition in \eqref{spatial_average_decay}).
Note that in the first inequality we used $(1+s)^{-2(1+\rho)(1-\zeta)+1}\leq (1+s)^{-(1+\rho)}$ for any $s>0$, as a result of the choice of $\zeta$ in \eqref{constraint}.

On the one hand, it follows from Lemma~\ref{longtime_separation}
and \eqref{eq:ls_inequality} that, for any $s\in A_M(T_M)$,
\begin{align*}
    & \AbsNormal{ E_{\mathrm{free}}(\phi(s),\psi(s)) - E_\infty }^{2(1-\zeta)}
    \chi_{A_M(T_M)}(s) \\
    & \quad \le \AbsNormal{ E_{\mathrm{free}}(\phi(s),\psi(s)) - E_\infty }^{2(1-\theta_{i(s)})} \AbsNormal{ E_{\mathrm{free}}(\phi(s),\psi(s)) - E_\infty }^{2(\theta_{i(s)}-\zeta)} \\
    & \quad \le C \BracketBig{ \NormNormal{\nabla \mu_\phi(s)}^2
    + \NormNormal{\nabla \mu_\psi(s)}^2
    + \AbsBig{\overline{\phi}(s)-\lim\limits_{q \to \infty}\overline{\phi(q)}}^{2(1-\theta_{i(s)})} } \\
    & \quad \le
    C \BracketBig{
    \NormNormal{\nabla \mu_\phi(s)}^2
    + \NormNormal{\nabla \mu_\psi(s)}^2
    + \AbsBig{\overline{\phi}(s)-\lim\limits_{q \to \infty}\overline{\phi(q)}} } \\
    & \quad \le
    C \BracketBig{
    \NormNormal{\nabla \mu_\phi(s)}^2
    + \NormNormal{\nabla \mu_\psi(s)}^2
    + (1+s)^{-2(1-\zeta)(1+\rho)} },
\end{align*}
where the index $i(s) \in \{1,\dots,n\}$ is chosen such that $(\phi(s),\psi(s)) \in B_{V \times V} ((\phi_{i(s)},\psi_{i(s)}),\delta)$, and we have used the fact that the residual energy term $\AbsNormal{ E_{\mathrm{free}}(\phi(s),\psi(s)) - E_\infty }^{2(\theta_{i(s)}-\zeta)}$ is uniformly bounded since $\zeta \le \theta_{i(s)}$. The last two inequalities follow from the fact that $2(1-\theta_{i(s)}) \ge 1$ and the first condition in \eqref{spatial_average_decay}, respectively. Indeed, it follows from \eqref{eq:spatial_average} (or equivalently, from the equation \eqref{eq:chdf3}) that
\begin{align*}
    \AbsNormal{\overline{\phi(q)} - \overline{\phi(s)}} & = \AbsMax{ \int_s^q \overline{\sigma_1(\phi)} \BracketBig{\overline{\phi}-c} \ \mathrm{d} \tau }
    \le \int_s^q \overline{\sigma_1(\phi)} \AbsBig{\overline{\phi}-c} \ \mathrm{d} \tau,
\end{align*}
for all $s,q \in [0,\infty)$.
Taking $q \to \infty$, we obtain
\[
\AbsNormal{\overline{\phi(s)} - \lim\limits_{q \to \infty} \overline{\phi(q)}} \le \int_s^\infty \overline{\sigma_1(\phi)} \AbsBig{\overline{\phi}-c} \ \mathrm{d} \tau \le C (1+s)^{-2(1+\rho)}.
\]
On the other hand, by the definition of the good times $A_M(T_M)$, we have that, in the set of ``bad times'', for almost any $s\in(t_*,\infty)$,
\begin{align*}
    & \AbsNormal{ E_{\mathrm{free}}(\phi(s),\psi(s)) - E_\infty }^{2(1-\zeta)}
    \chi_{(t_*,\infty) \backslash A_M(T_M)}(s)
    \le
    C \frac{
    \NormNormal{\nabla \mu_\phi(s)}^2
    + \NormNormal{\nabla \mu_\psi(s)}^2}{M^2}.
\end{align*}
Consequently, it holds that
\begin{align*}
    & \BracketMax{
    \int_s^\infty
    \NormNormal{\nabla \mu_\phi}^2
    + \NormNormal{\nabla \mu_\psi}^2
    + \NormNormal{\sp}_{\bm{L}^2_\sigma(\Omega)}^2
    + \NormNormal{\sp}_{\bm{L}^r_\sigma(\Omega)}^r
    +(1+\tau)^{-2(1-\zeta)(1+\rho)} \ \mathrm{d} \tau
    }^{2(1-\zeta)} \\
    & \quad \le
    C \BracketBig{
    \NormNormal{\nabla \mu_\phi(s)}^2
    + \NormNormal{\nabla \mu_\psi(s)}^2
    + \NormNormal{\sp(s)}_{\bm{L}^2_\sigma(\Omega)}^2
    + \NormNormal{\sp(s)}_{\bm{L}^r_\sigma(\Omega)}^r
    + (1+s)^{-2(1-\zeta)(1+\rho)} },
\end{align*}
for almost any $s\in(t_*,\infty)$.
We then apply \cite[Lemma 7.1]{FS} (see also \cite[Lemma~A.1]{GP2025}) to conclude that
\begin{equation} \label{L1_integrability}
    |\nabla \mu_\phi|, |\nabla \mu_\psi|, \AbsNormal{\sp}, \AbsNormal{\sp}^{\frac{r}{2}}
    \in L^1\BracketBig{t_*,\infty;L^2(\Omega)}.
\end{equation}
Thus, a comparison argument in equations \eqref{eq:chdf3} and \eqref{eq:chdf5} yields
\[
\partial_t \phi, \partial_t \psi
\in L^1 \BracketBig{t_*,\infty;V^*}.
\]
 The remainder of the proof proceeds as in \cite{GP2025}. The proof is finished.
\qed

\section*{Declarations}
\noindent
\textbf{Conflict of interest.} The authors have no competing interests to declare that are relevant to the content of this article.
\smallskip
\\
\textbf{Fundings.}
The research of AP was funded in part by the Austrian Science Fund 
\href{https://doi.org/10.55776/ESP552}{10.55776/ESP552}.
The research of HW was partially supported by the Natural Science Foundation of Shanghai (No. 25ZR1401023).
\smallskip
\\
\noindent
\textbf{Data availability.} Data sharing not applicable to this article as no datasets were generated or analyzed during the current study.
\smallskip
\\
\noindent
\textbf{Acknowledgments.}
Part of these results have been obtained as a consequence of some stimulating discussions between MG, AP, and HW during the Thematic Program on \textit{Free Boundary Problems} at the Erwin Schrödinger International Institute for Mathematics and Physics (ESI) in Vienna, whose hospitality is kindly acknowledged.
MG and AP are members of Gruppo Nazionale per l'Ana\-li\-si Matematica, la Probabilit\`{a} e le loro Applicazioni (GNAMPA), Istituto Nazionale di Alta Matematica (INdAM). MG's research is part of the activities of ``Dipartimento di Eccellenza 2023--2027'' of Politecnico di Milano.
Ouyang is supported by the China Scholarship Council (CSC) Program (Grant No. 202506100096).
HW is a member of Key Laboratory of Mathematics for Nonlinear Sciences (Fudan University), Ministry of Education of China.
For open access purposes, the authors have applied a CC BY public copyright license to
any author accepted manuscript version arising from this submission.

\end{document}